\documentclass[11pt]{amsart}
\usepackage{stmaryrd}
\usepackage{dynkin-diagrams}
\usepackage{tikz-cd}
\usepackage{etoolbox}
\usepackage[flushleft]{threeparttable}
\usepackage{hhline}
\usepackage{comment}
\usepackage{tikz}
\usepackage[page]{appendix}
\usepackage[utf8]{inputenc}
\usepackage{filecontents}
\usepackage{graphicx} 
\usepackage{mathtools}
\usepackage{amsmath}
\usepackage{amsthm}
\usepackage{amsfonts}
\usepackage{amssymb}
\usepackage{mathrsfs}
\usepackage{xypic}
\usepackage{color}
\usepackage{xcolor}
\usepackage{titlesec}
\usepackage{titletoc}
\usepackage{tabularx}
\usepackage[colorlinks]{hyperref}
\newtheorem{lem}{Lemma}
\newtheorem{conj}{Conjecture}
\newtheorem{cor}{Corollary}
\newtheorem{prop}{Proposition}
\newtheorem{thm}{Theorem}
\newtheorem*{quest*}{Question}
\newtheorem*{thm*}{Theorem}
\newtheorem{defn}{Definition}

\newcommand{\ev}{\operatorname{ev}}
\newcommand{\tf}{\operatorname{tf}}
\newcommand{\pcrys}{\operatorname{pcrys}}
\newcommand{\Frob}{\operatorname{Frob}}

\usepackage{bbm}
\usepackage[colorlinks]{hyperref}
\definecolor{cclr}{rgb}{25,25,112}
\hypersetup{
citecolor=[rgb]{0.15, 0.15, 0.68},
linkcolor=gray,   
urlcolor=magenta}

\newcommand{\Ch}{\operatorname{Ch}}
\newcommand{\SP}{\operatorname{sp}}
\newcommand{\Top}{{\operatorname{top}}}

\newcommand{\Sha}{{\operatorname{Sha}}}
\newcommand{\fl}{\operatorname{fl}}
\newcommand{\wt}[1]{\widetilde{#1}}
\newcommand{\sS}{\operatorname{ss}}
\newcommand{\BKF}{\operatorname{BKF}}

\newcommand{\frep}{^f\operatorname{Rep}}

\newcommand{\id}{\operatorname{id}}
\newcommand{\Ga}{\operatorname{G}_a}
\newcommand{\cC}{\mathcal{C}}

\newcommand{\fM}{\mathfrak{M}}

\newcommand{\clR}{\mathcal{R}}

\newcommand{\Inf}{{\operatorname{inf}}}
\newcommand{\fS}{\mathfrak{S}}
\newcommand{\st}{{\operatorname{st}}}

\newcommand{\Aut}{\operatorname{Aut}}
\newcommand{\red}{\operatorname{red}}

\newcommand{\IsomFil}{\operatorname{Isom-fil}}
\newcommand{\wh}{\breve}

\newcommand{\can}{\operatorname{can}}

\newcommand{\Vect}{\operatorname{Vect}}

\newcommand{\gr}{\operatorname{gr}}
\newcommand{\Res}{\operatorname{Res}}

\newcommand{\Rep}{\operatorname{Rep}}

\newcommand{\spec}{\operatorname{Spec}}
\newcommand{\spf}{\operatorname{Spf}}

\newcommand{\dirlim}[1]{\underset{#1}{\underrightarrow{\operatorname{lim}}}}
\newcommand{\invlim}[1]{\underset{#1}{\underleftarrow{\operatorname{lim}}}}
\usepackage{xstring}

\titleformat{\section}[runin]{\normalfont\bfseries}{\thesection.}{3pt}{}
\titleformat{\subsection}[runin]{\normalfont\bfseries}{\thesubsection.}{3pt}{}
\titleformat{\subsubsection}[runin]{\normalfont\bfseries}{\thesubsubsection.}{3pt}{}
\renewcommand{\thesection}{\arabic{section}}
\titleformat{\section}{\normalfont\large\bfseries}{\thesection.~~}{1em}{}
\DeclareMathAlphabet{\mathpzc}{OT1}{pzc}{m}{it}

\dottedcontents{section}[2.5em]{}{2.9em}{1pc}

\begin{document}

\newcommand{\tor}{\text{tor}}
\newcommand{\cf}{\text{cf}}
\newcommand{\chicyc}{\chi_{\operatorname{cyc}}}
\newcommand{\Brnr}{\operatorname{Br}_{\operatorname{nr}}}
\newcommand{\Bronr}{\operatorname{Br}_{1,\operatorname{nr}}}
\newcommand{\Br}{\operatorname{Br}}
\newcommand{\LT}{\operatorname{LT}}
\newcommand{\Qp}{\mathbb{Q}_p}
\newcommand{\TODO}{{\color{red} TODO}}
\newcommand{\Gm}{\mathbb{G}_m}
\newcommand{\Scin}{\operatorname{Scin}}
\newcommand{\Fil}{\operatorname{Fil}}
\newcommand{\Ind}{\operatorname{Ind}}
\newcommand{\Sym}{\operatorname{Sym}}
\newcommand{\sym}{\operatorname{sym}}
\newcommand{\semis}[1]{#1\operatorname{-ss}}
\newcommand{\sem}{\operatorname{ss}}
\newcommand{\Alt}{\operatorname{Alt}}
\newcommand{\DGK}{\operatorname{DG}_K}
\newcommand{\Img}{\operatorname{Im}}
\newcommand{\Ker}{\operatorname{Ker}}
\newcommand{\grp}{\operatorname{grp}}
\newcommand{\cont}{\operatorname{cont}}
\newcommand{\ord}{\operatorname{ord}}
\newcommand{\Hom}{\operatorname{Hom}}
\newcommand{\Ext}{\operatorname{Ext}}
\newcommand{\Ad}{\operatorname{Ad}}
\newcommand{\Id}{\operatorname{id}}
\newcommand{\Lie}{\operatorname{Lie}}
\newcommand{\Lift}{\operatorname{Lift}}
\newcommand{\ad}{\operatorname{ad}}
\newcommand{\rk}{\operatorname{rk}}
\newcommand{\Det}{\operatorname{Det}}
\newcommand{\LHS}{\operatorname{LHS}}
\newcommand{\RHS}{\operatorname{RHS}}
\newcommand{\bFp}{\bar{\mathbb{F}}_p}
\newcommand{\bZp}{\bar{\mathbb{Z}}_p}
\newcommand{\bQp}{\bar{\mathbb{Q}}_p}
\newcommand{\bF}{{\mathbb{F}}}
\newcommand{\bZ}{{\mathbb{Z}}}
\newcommand{\bA}{{\mathbb{A}}}
\newcommand{\bB}{{\mathbb{B}}}
\newcommand{\bP}{{\breve{P}}}
\newcommand{\WD}{\operatorname{WD}}
\newcommand{\bM}{\breve{M}}
\newcommand{\bT}{\breve{T}}
\newcommand{\bS}{\breve{S}}
\newcommand{\brB}{\breve{B}}
\newcommand{\bG}{\breve{G}}
\newcommand{\bH}{\breve{H}}
\newcommand{\bQ}{{\mathbb{Q}}}
\newcommand{\bR}{{\mathbb{R}}}
\newcommand{\bC}{{\mathbb{C}}}
\newcommand{\scrC}{{\mathscr{C}}}
\newcommand{\Fq}{\bar{\mathbb{F}}_q}
\newcommand{\GSp}{\operatorname{GSp}}
\newcommand{\Sp}{\operatorname{Sp}}
\newcommand{\GL}{\operatorname{GL}}
\newcommand{\cris}{\operatorname{cris}}
\newcommand{\modp}{\operatorname{mod} \varpi}
\newcommand{\End}{\operatorname{End}}
\newcommand{\SL}{\operatorname{SL}}
\newcommand{\SO}{\operatorname{SO}}
\newcommand{\GO}{\operatorname{GO}}
\newcommand{\cO}{\mathcal{O}}
\newcommand{\cX}{\mathcal{X}}
\newcommand{\cY}{\mathcal{Y}}
\newcommand{\cN}{\mathcal{N}}
\newcommand{\cU}{\mathcal{U}}
\newcommand{\cV}{\mathcal{V}}
\newcommand{\cZ}{\mathcal{Z}}
\newcommand{\Lieg}{\mathfrak{g}}
\newcommand{\Lieh}{\mathfrak{h}}
\newcommand{\Mat}{\operatorname{Mat}}
\newcommand{\Gal}{\operatorname{Gal}}
\newcommand{\Art}{\operatorname{Art}}
\newcommand{\BHT}{\mathbb{B}_{\operatorname{HT}}}
\newcommand{\DXT}{{\hatI}^{X(T)}}
\newcommand{\HTGC}{\prod_{\sigma:K_m\hookrightarrow \bC}X_*(G_{\bC})}
\newcommand{\HTGCa}{\prod_{\tilde\tau\in S}X_*(G_{\bC})}
\newcommand{\HT}{{\operatorname{HT}}}
\newcommand{\hT}{{\mathcal{HT}}}
\newcommand{\dR}{{\operatorname{dR}}}
\newcommand{\hatbI}{\widehat{\bar I}_K}
\newcommand{\hatI}{\widehat{I}_K}
\newcommand{\symt}{\operatorname{sym}^2}
\newcommand{\altt}{\operatorname{alt}^2}
\newcommand{\socle}{\operatorname{soc}}
\newcommand\nlcup{%
  \mathrel{\ooalign{\hss$\cup$\hss\cr%
  \kern0.7ex\raise0.6ex\hbox{\scalebox{0.4}{$\diamond$}}}}}
\newcommand{\matt}[9]{
\left(
\begin{matrix}
#1 & #2 & #3 \\
#4 & #5 & #6 \\
#7 & #8 & #9
\end{matrix}
\right)
}
\newcommand{\lsup}[2]{{}^{#1}\hspace{-0.6mm}#2}

\author{Lin, Zhongyipan}
\title{Extensions of de Rham Galois representations}

\begin{abstract}
We construct the parabolic version and the reductive
version of 
the integral
de Rham moduli stacks of Langlands parameters
($p>3$).
We allow the group to be arbitrarily ramified.

We propose that the top Chow group of the reduced
Emerton-Gee stack $\cX_{\lsup LG}$
is isomorphic to
that of the moduli
of Weil-Deligne representations valued in $\lsup LB$,
where $\lsup LB$ is a Borel of $\lsup LG$.
The latter bears a concrete description
by Serre weights corrected by the Kottwitz homomorphism.

We explicitly define such a map using parabolic de Rham
moduli stacks
as the composition of 
a chain of tautological maps,
and confirm it is an isomorphism
for
(1) algebraic tori, (2) unitary, orthogonal and symplectic
groups,
(3) tame groups when restricted to
the cyclotomic-free part of the Chow groups.
\end{abstract}

\maketitle

Let $K/\Qp$ be a local field with absolute Galois group $\Gal_K$.
Let $(G, B, T, \{X_\alpha\})$ be a pinned quasi-split
connected reductive group over $K$ that splits over $L$,
with dual pinned group
$(\bG, \brB, \bT, \{Y_\alpha\})$
and Langlands dual
$\lsup LG = \bG \rtimes \Gamma$
where $\Gamma=\Gal(L/K)$.
Denote by $M$ a standard Levi subgroup of $G$.
Both $M$ and its dual $\bM$
have pinnings induced from the pinning
of $G$ and $\bG$, respectively.
The pinned homomorphism
$\bM\to \bG$ induces
$\lsup LM =\bM\rtimes \Gamma\hookrightarrow\lsup LG$,
and we choose the parabolic $P=U\rtimes \bM$
of $\bG$ by demanding $\brB\subset P$.
Set $\lsup LP:=P\rtimes \Gamma$.

\begin{conj}
\label{conj}
There is a natural isomorphism
between the top Chow group
of the reduced Emerton-Gee stacks
$\cX_{K, \lsup LG, \red}$
and
the top Chow group
of the moduli
of Weil-Deligne representations
$\WD_{K, \lsup LB}$
valued in $\lsup LB$.
\end{conj}

\begin{thm} [Proposition \ref{prop:main}]
\label{thm:main}
If $p>2$,
the Conjecture holds
when
(1) $G$ is an algebraic torus,
(2) $G=U_n$ is a unitary group, or
(3) $\lsup LG=\GSp_{2n}$
or $\GO_n$.
\end{thm}

In general, it is harmless
to replace $G$ by $\Res_{K/\bQ_p}G$
and assume $K=\bQ_p$;
we will suppress the appearance of $K$
in notations when it is clear.
Suppose $G$ admits
a local twisting element $-\rho$
(in the sense of \cite{GHS}).
We consider the moduli
of de Rham parameters $\cX_{\lsup LB}^{\rho}$
over $\spf \bar\bZ_p$.
The expectation is that
all maps in the following
diagram
$$
\xymatrix{
\cX_{\lsup LG, \red}
 &
\cX_{\lsup LB, \red}^{\rho}
\ar@{^{(}->}[r]^i
\ar[l]_{f}
 &
\cX_{\lsup LB}^{\rho}
 & 
 \cX^{\rho}_{\lsup LB}[\frac{1}{p}]
\ar@{_{(}->}[l]_j
\ar[r]^{\WD}
&
\WD_{\lsup LB}
}
$$
induce an isomorphism
of top Chow groups,
and thus
$\Ch_{\Top}(\WD_{\lsup LB})
\xrightarrow{\SP_f\circ i^*\circ (j^*)^{-1}\circ\SP_{\WD}^{-1}}
\Ch_{\Top}(\cX_{\lsup LG, \red})$
is an isomorphism.
Here $\SP_{f}([C]):=[f(C)]$,
$i^*$ is the usual Gysin homomorphism
(\cite[Tag 02T7]{Stacks})
and $j^*$ is the usual flat pullback.

\begin{thm}
[Definition \ref{def:XP}, Theorem \ref{thm:XP}]
\label{thm:moduli}
Let $\lambda$ be any Hodge type $\Gm\to \bT$
that pairs negatively with roots of $U$
and let $\tau$ be any inertial type.
If $p>3$, then
$\cX_{\lsup LP}^{\tau,\lambda}$
is a $\bar\bZ_p$-flat
Noetherian formal algebraic stack
over $\spf \bar\bZ_p$.
If $\lsup LP=\lsup LG$,
then
$\cX_{\lsup LG}^{\tau,\lambda}$
is of finite type over $\spf \bar\bZ_p$.
\end{thm}

\subsection{Remark}
(1)
``Pairs negatively'' is the condition that allows
$\cX_{\lsup LP}^{\tau, \lambda}$
to achieve maximal dimension.
Theorem 
\ref{thm:moduli}
as well as almost all of its supporting lemmas
and propositions
still holds
after replacing ``pairs negatively''
by ``pairs non-vanishingly''
(some roots positive, some roots negative,
but none zero)
with minimal extra effort.
But we refrain from introducing the complexity
as the more general cases are not very useful to us.

(2)
There is evidence suggesting that
$\cX_{\lsup LP}^{\tau, \lambda}$
is of finite type over $\spf \bar\bZ_p$
in general, but we don't need this.

(3)
It is the first work in which
de Rham moduli stacks
for reductive groups
and de Rham deformation rings
for parabolic groups
are constructed
(versal rings of $\cX_{\lsup LP}^{\tau, \lambda}$
furnish such deformation rings),
in generality.

\subsection{}
The moduli stack $\WD_{\lsup LB}$
admits very concrete description,
and we collect some of its properties
that are clean to formulate.

\begin{thm}[Proposition \ref{prop:sp},
Lemma \ref{lem:comb},
Lemma \ref{lem:kottwitz}]
\label{thm:sp}~

(1)
$\SP_{\WD}$ is an isomorphism
of top Chow groups.

(2)
Let $\Delta=\Delta_G(B)=\{\alpha_1, \dots, \alpha_s\}$
be the set of simple positive roots.
Then
$\WD_{\lsup LB}\cong
\WD_{\prod_{i=1}^sU_{\sigma_i}\rtimes \lsup LT}$.

(3)
The trivial monodromy part of
$\WD_{\lsup LB}$
is isomorphic to
$\WD_{\lsup LT}$.

(4)
Write $\Ch_{\Top}(\WD_{\lsup LT})=\bZ^{\oplus I}$.
Then the index set $I$ has a torsion abelian group structure
fitting in
the short exact sequence
$$
0\to
\Hom(X^*(T)_{I_{\bQ}, \tor}^\Gamma, \bQ/\bZ)
\to
I
\to
\Hom(X^*(T)_{I_{\bQ}, \tf}^\Gamma, (\bQ/\bZ)_{p'})
\to 0,
$$
which gives a combinatorial description
of potentially crystalline components
of $\cX_{\lsup LB}^{\rho}$.
\end{thm}

We remark that 
$\Hom(X^*(T)_{I_{\bQ}, \tf}^\Gamma, (\bQ/\bZ)_{p'})$
is the set of regular Serre weights,
and
$I^\vee \to X^*(T)_{I_{\bQ}, \tor}^\Gamma$
is the Kottwitz homomorphism.

We are currently working on a
computer-assisted
proof that resolves
Conjecture \ref{conj}
in full generality, at least for large $p$.
This project takes a lot of coding effort and is not
expected to be finished anytime soon.
In this note, we make minimal efforts
to establish the cyclotomic-free case
as it is used by some other projects.

\begin{thm}[Proposition \ref{prop:LG-LT}]
Assume $G$ is tame.
Denote by $\cX^\cf_{\lsup LG}$
the substack
of $\cX_{\lsup LG}$
defined by the condition $\{\rho|H^2(\Gal_K, \ad\rho)=0\}$.
Then the chain
$
\cX^{\cf}_{\lsup LG, \red}
\leftarrow
\cX^{\cf}_{\lsup LB, \red}
\to
\cX^{\cf}_{\lsup LT, \red}
$
induces isomorphism of top Chow groups.
\end{thm}

Finally, we also present the following
technical result that has independent interest.

\begin{thm} [Corollary \ref{cor:semistable}, Corollary \ref{cor:de-rham}]
Let $\rho: \Gal_{K}\to \lsup LM(\bQp)$ be a de Rham representation whose Hodge type is $P^-$-dominant.
Endow $\Lie U(\bQp)$ with $\Gal_K$-action 
via $\ad\rho$.

Any $\rho':\Gal_K\to \lsup LP(\bQp)$
such that $\rho'\times^P\bM= \rho$
is also de Rham
with its possible Weil-Deligne
types parameterized by $H^2(\Gal_K, U(\bQp))$.
\end{thm}

\subsection{Acknowledgement}
The author thanks Bao Le Hung and Tony Feng
for asking for a sane description
of the top Chow group of reduced Emerton-Gee stacks,
which is the source of this work.

\tableofcontents

\section{Preliminaries}~

For ease of notation, we write $-\otimes-$ for
$-\otimes_{\Qp}-$.

\subsection{Exact $\otimes$-filtrations}
Let $\Vect_A$ be the category of projective $A$-modules.
For $X\in\Vect_A$, a filtration of $X$ indexed by $\bZ$
is a tuple $(\Fil^nX)_{n\in \bZ}$ where $\Fil^nX\in\Vect_A$, $\Fil^nX \supset \Fil^{n+1}X$,
$\cap \Fil^n X = 0$, and $\cup \Fil^nX=X$.
Denote by $\Fil\Vect_A$ the category of filtered (indexed by $\bZ$) projective $A$-modules and denote by $\gr\Vect_A$ be the category of $\bZ$-graded projective $A$-modules.

Let $\scrC$ be an ind-tannakian category (\cite[III 1.1.1]{SN72}) over a ring $A$.
Let $\omega: \scrC\to\Vect_A$ be an exact tensor functor.
An \textit{exact $\otimes$-filtration} $\theta$ on $\omega$ is a factorization
$$
\xymatrix{
    \scrC \ar[rr]^{\omega}\ar[rd]^{\theta} & & \Vect_A \\
    & \Fil\Vect_A \ar[ru]^{\text{forget}} &
}
$$
such that the following are satisfied:
\begin{itemize}
\item[(FE 1)]
For $X\in \scrC$, $\theta^n(X)$ is a direct summand of $\omega(X)$;
\item[(FE 2)]
The associated graded functor $\gr(\theta)$ is exact.
\item[(FE 3)]
For all $n\in\bZ$, $X, Y\in\scrC$,
$$
\theta^n(X\otimes Y) = \sum_{i+j=n}\theta^i(X)\otimes\theta^j(Y).
$$
\end{itemize}

Let $\omega$, $\omega'$ be exact tensor functors with
exact $\otimes$-filtrations $\theta$, $\theta'$, respectively.
Denote by
$
\IsomFil^{\otimes}(\theta, \theta')
$
the set of tensor isomorphisms $f:\omega \to \omega'$ inducing an isomorphism of filtrations $f_*\theta \cong \theta'$.
Set
$$
\Aut^{\otimes}(\theta) := \IsomFil^{\otimes}(\theta, \theta)
$$
and denote by
$\Aut^{\otimes!}(\theta)$
the subfunctor of
$\Aut^{\otimes}(\theta)$ which induces identity of the associated grading.

\vspace{3mm}

\subsection{Splitting of exact $\otimes$-filtrations}

An exact tensor functor
$\lambda: \scrC\to\gr\Vect_A$
induces a canonical exact $\otimes$-filtration $\lambda_{\can}$, which is defined as
$\lambda^n_{\can}(X):= \sum_{n' \ge n}\lambda^{n'}(X)$ for all $X\in\scrC$.

An exact $\otimes$-filtration is
said to be \textit{splittable}
if it is isomorphic to the
    canonical exact $\otimes$-filtration
    associated to a graded exact tensor functor.

\subsection{Torsors for smooth affine groups}
For smooth affine groups $H$,
we interpret $H$-torsors over $\spec A$
as rigid exact tensor functors $\Rep(H)\to \Vect_A$.
We specialize to the case where $A=K\otimes\bC$.
Here are some basic facts
\begin{enumerate}
\item 
Any exact $\otimes$-filtration on
$H$-torsors over $\spec A$ is splittable
(\cite[IV.2.4, IV.1.3]{SN72}).
\item
Any exact $\otimes$-grading on
$H$-torsors over $\spec A$ is the weight grading
of a cocharacter $\lambda: (\Gm)_A \to H_A$
(\cite[IV.1.3]{SN72}).
\item
As a consequence of (1) and (2),
any exact $\otimes$-filtration on
$H$-torsors over $\spec A$ is equal to
$\lambda_{\can}$
for some cocharacter $\lambda: (\Gm)_A \to H_A$.
\end{enumerate}
Let $\omega$ be an $H$-torsor.
For each group scheme morphism $f:H\to H'$,
we can functorially define 
the pushforward $f_*\omega=\omega \times^H H'$,
which is an $H'$-torsor.
The construction of pushforward extends to exact $\otimes$-filtrations and $\otimes$-gradings:
if $\lambda:(\Gm)_A\to H_A$ is an exact $\otimes$-grading, then $f_*(\lambda)=f \circ \lambda$,
and if $\theta=\lambda_{\can}$ is an exact $\otimes$-filtration, then
$f_*(\theta)=(f_*\lambda)_{\can}$.

In particular, if $f=\Ad(g)$, $g\in H(A)$
is an inner automorphism,
then $f_*(\lambda_{\can})=(g\lambda g^{-1})_{\can}$,
which we also denote by $g(\lambda_{\can})g^{-1}$
by abuse of notation.

\begin{lem}
\label{triv-fil}
If $\theta$ is an exact $\otimes$-filtration
on an $P$-torsor over $K\otimes\bC$,
then there exists an element 
$u\in U$
such that
$$
u\theta u^{-1} = (\theta \times^P\bM) \times^{\bM}P.
$$
\end{lem}

\begin{proof}
Set $A=K\otimes \bC$.
Write $\theta=\lambda_{\can}$, where $\lambda:(\Gm)_{A}\to P_A$
is a cocharacter.
The image of $\lambda$ is contained in a maximal torus $S$ of $P_A$.
There exists an element $u\in U(A)$
such that $u^{-1} S u \subset \bM_A$.
Replace $\theta$ by $u \theta u^{-1}$,
so $S\subset \bM_A$ and
$\lambda = (\lambda \times^P\bM) \times^{\bM}P$.
\end{proof}

\subsection{Dynamic methods}
\label{dyn-fil}
Since $F \mapsto F(K\otimes \bC)$
is a fully faithful functor
from the category of scheme-theoretic torsors over $K\otimes \bC$ to the category of set-theoretic torsors,
we will not distinguish
$F$ and $F(K\otimes \bC)$
in the entire article.

Recall some notations from \cite[Section 4.1]{Crd11}.
Let $\lambda$ be a cocharacter of a reductive group $H$ over $K\otimes \bC$.
Define
$P_H(\lambda):=\{g\in H|\lim_{t\to 0}\lambda(t)g\lambda(t)^{-1}\text{~exists.}\}$;
$P_H(\lambda)$ is a smooth subgroup of $H$,
and all parabolic subgroups of $H$ are of the form
$P_H(\lambda)$ for some $\lambda$.
Define $U_H(\lambda)=\{g\in H|\lim_{t\to 0}\lambda(t)g\lambda(t)^{-1}=1\}$. Then $U_H(\lambda)\subset P_H(\lambda)$  is the unipotent radical.

\vspace{3mm}
Recall that $\lambda_{\can}=(\lambda_{\can}^i)_{i\in \bZ}$ is the exact $\otimes$-filtration associated to $\lambda$.

\vspace{3mm}

\begin{thm}
Consider the adjoint representation $\ad: H\to \GL(\Lie(H))$.
We have
\begin{align*}
\Lie \Aut^{\otimes}(\lambda_{\can})
&= \lambda^0_{\can}(\Lie(H)), ~\text{and}\\
\Lie \Aut^{\otimes!}(\lambda_{\can})
&= \lambda^1_{\can}(\Lie(H)).
\end{align*}

As a consequence, 
we have $P_H(\lambda) = \Aut^{\otimes}(\lambda_{\can})$
and $U_H(\lambda) = \Aut^{\otimes!}(\lambda_{\can})$.
\end{thm}

\begin{proof}
The first paragraph is a special case of \cite[IV 2.1.4.1]{SN72} where $\alpha=0,1$.
The second paragraph follows from \cite[Theorem 4.1.7(4)]{Crd11}.
\end{proof}

\section{Extension of weakly admissible filtered $(\phi, N)$-modules}

In this section, we work over a general field $K$
(instead of specializing to $\bQ_p$)
and we always work with connected groups.

\begin{defn}
Let $K/\Qp$ be a finite extension.
Set $K_0 = W(k)[1/p]$ where $k$ is the residue field of $K$.
A \textit{filtered $(\phi, N)$-module with $\bC$-coefficients} is a tuple
$(D, \phi_D, \theta_D, N_D)$ where
\begin{itemize}
\item[-]
$D$ is a finite free module over 
$K_0\otimes \bC$;
\item[-]
$\phi_D: (\phi\otimes 1)^*D\to D$ is an isomorphism of $K_0\otimes \bC$-modules;
\item[-]
$\theta_D$ is a filtration on $D_K:=D\otimes_{K_0}K$ such that
$\theta_D^jD_K=0$ if $j\gg 0$, and $\theta_D^jD_K=D_K$ if $j\ll 0$.
\item[-]
$N_D:D\to D$ is a $K_0$-linear
endomorphism
such such $N_D\phi_D=p\phi_DN_D$.
\end{itemize}
Here $\phi\otimes 1:K_0\otimes \bC\to K_0\otimes \bC$ sends $x\otimes y$ to $\varphi(x)\otimes y$.
\end{defn}

We will drop the phrase ``with $\bC$-coefficients''
because it is always enforced.

\begin{defn}
Let $H$ be an affine smooth group over $\bC$.
A \textit{filtered $(\phi, N)$-module and $H$-structure} is a tuple
$(F, \phi_F, \theta_F, N_F)$ such that
\begin{itemize}
\item[-]
$F$ is a $H$-torsor over $\spec K_0\otimes \bC$;
\item[-]
$\phi_F: (\phi\otimes 1)^*F\to F$ is a $H$-equivariant isomorphism over $\spec K_0\otimes \bC$;
\item[-]
$\theta_F$ is an exact $\otimes$-filtration 
on $F\otimes_{K_0}K$;
\item[-]
$N_F\in \Lie(H)\otimes_{\bQ_p}K_0$
such that
$N_F = p \phi_F (\phi\otimes 1)^{*} N_F \phi_F^{-1}$.
\end{itemize}
\end{defn}

We will drop the phrase ``with $H$-structure''
when the choice of $H$ is clear from the context.

\subsection{Weak admissibility and semistable representations}
A filtered $(\phi, N)$-module $(F, \phi_F, \theta_F, N_F)$
with $H$-structure 
is said to be \textit{weakly admissible}
if for any algebraic representation $H\to \GL(V)$, 
$(F, \phi_F, \theta_F, N_F)\times^HV$ is 
a weakly admissible filtered $(\phi, N)$-module.

Since the covariant Fontaine's functors $V_{\operatorname{ss}}$ and $D_{\operatorname{ss}}$ are rigid exact tensor functors (see the paragraph before \cite[9.1.9]{C11}),
the category of weakly admissible filtered $(\phi, N)$-module with $H$-structure
is equivalent to the category of semistable representations valued in $H$.

\begin{lem}
(1) If $(F, \phi_F, \theta_F, N_F)$ is weakly admissible,
then $(F, \phi_F, \theta_F, N_F)\times^HH'$
is also weakly admissible for any group
scheme morphism $H\to H'$.

(2) If $H\to H'$ is an embedding, then
the converse is also true.
\label{lemma:aw-pushforward}
\end{lem}

\begin{proof}
(1) Unravelling definitions shows $(F, \phi_F, \theta_F, N_F)\times^HH'
= D_{\operatorname{ss}}(V_{\operatorname{ss}}((F, \phi_F, \theta_F, N_F)) \times^HH')$, which is automatically weakly admissible.

(2) 
It is harmless to replace $H'$ by $\GL(V)$.
It follows from the fact that a single faithful
representation of $H$ generates
the Tannkian category of algebraic representations $\Rep(H)$.
\end{proof}

Recall that $P$ is a parabolic subgroup of $\bG$
with Levi subgroup $\bM$ and unipotent radical $U$.

\begin{prop}
Let $(F, \phi_F, \theta_F, N_F)$ be a filtered $(\phi, N)$-module
with $P$-structure.
Then
$(F, \phi_F, \theta_F, N_F)$ is weakly admissible
if and only if
$(F, \phi_F, \theta_F, N_F)\times^P\bM\times^{\bM}P$
is weakly admissible.
\label{prop:weak-adm-P}
\end{prop}

\begin{proof}
Since $P\to \bM \times P/Z(\bG)$
is an embedding,
Lemma \ref{lemma:aw-pushforward} allows us to replace $P$, $\bM$,
and $\bG$ by $P/Z(\bG)$, $\bM/Z(\bG)$ and $\bG/Z(\bG)$.
So now $P\subset \bG\to \GL(\Lie \bG)$ is an embedding.
For ease of notation, set $(-)^{\semis{\bM}}:=(-)\times^P\bM\times^{\bM}P$,
and set $F_*:=(F, \phi_F, \theta_F, N_F)$.

Choose a cocharacter $\mu$ of $P$ such that
$P=P_{\bG}(\mu)$.
Then $P':=P_{\GL(\Lie \bG)}(\mu)\supset P$
is a parabolic subgroup of $\GL(\Lie \bG)$
with unipotent radical $U'=U_{\GL(\Lie \bG)}(\mu)\supset U$
and Levi subgroup $\bM':=Z_{\GL(\Lie\bG)}(\mu)\supset \bM$.
Since
$$
\xymatrix{
    U \ar[d]\ar[r] & U'\ar[d] \\
    P \ar[r] & P'
}
$$
is Cartesian, we have
$$
(F_*\times^PP')^{\semis{M'}}
=
F_*^{\semis{\bM}}\times^PP'.
$$
$F_*^{\semis{\bM}}\times^PP'$ is weakly admissible by
Lemma \ref{lemma:aw-pushforward},
and therefore $F_*\times^PP'$
is weakly admissible
by \cite[Proposition 8.2.10]{BC08}
(which is the special case of this lemma
for $\GL_N$).
By Lemma \ref{lemma:aw-pushforward} again,
we conclude $F_*$
is weakly admissible.
\end{proof}

\subsection{Newton polygon of isocrystals}
\label{newton-poly}
Let $\breve{K}$ be the $p$-adic completion of the maximal unramified extension of $K_0$.
By the Dieudonn\'e-Manin classification, the category of isocrystals over $\breve K$ is a semisimple category.
The simple objects can be classified by rational numbers $s/r$, where $r$ is a positive integer and $s$ is an integer coprime to $r$.
Denote by $D_{r,s}$ the simple object labeled by the rational number $s/r$. $D_{r,s}$ has dimension $r$,
and we call $s/r$ the \textit{slope} of $D_{r,s}$.

Let $(D, \phi)$ be an isocrystal over $K_0$.
Then $\breve{D} = \breve{K} \otimes_{K_0} D$
is a direct sum of simple objects $D_{r_i,s_i}$.
We call the numbers $s_i/r_i$ that appear in the direct sum decomposition the slopes of $D$.
Say $D$ has slopes $\{\alpha_0<\cdots<\alpha_n\}$ with multiplicities $\{\mu_0,\cdots,\mu_n\}$.
The \textit{Newton polygon} of $D$ is the convex polygon with leftmost endpoint $(0,0)$, and having $\mu_i$ consecutive segments of horizontal distance $1$
and slope $\alpha_i$.

\begin{lem}
If all slopes of $D$ are positive numbers, then for any lattice $\mathcal{L}\subset D$, we have
$$
\lim_{n\to\infty}\phi^n\mathcal{L} = \{0\}
$$
in the sense that the diameter of the bounded sets 
$\phi^n\mathcal{L}$ converges to $0$.
Note that $\mathcal{L}$ is not assumed to be $\phi$-stable.
\label{lem:newton-poly}
\end{lem}

\begin{proof}
Let $N$ be the product of the denominator of the slopes of $D$.
By the Diedonn\'e-Manin classification,
there is a basis $\{x_1,..,x_t\}$ of
$\breve{K}\otimes_{K_0}D$
such that $\phi^Nx_i = p^{S_i}x_i$ for positive integers $S_i$, $1\le i\le t$.
The $p$-adic topology on $D$ can be defined by
$$
\vert \lambda_1x_1 +\cdots+\lambda_tx_t\vert
=\max_{1\le i\le t}(\vert \lambda_i\vert)
$$
where $\lambda_i\in \breve K$, $1 \le i\le t$.
For any $x\in D$, $|\phi^Nx| < \frac{1}{p}|x|$.
Therefore for any lattice $\mathcal{L}$, we have
$
\lim_{n\to\infty}\phi^{nN}\mathcal{L} = \{0\}
$.
Replacing $\mathcal{L}$ by $\phi^{k}\mathcal{L}$, $1\le k < N$, we have
$
\lim_{n\to\infty}\phi^{k + nN}\mathcal{L} = \{0\}
$.
Combining these, we have
$
\lim_{n\to\infty}\phi^{n}\mathcal{L} = \{0\}
$.
\end{proof}

\begin{cor}
\label{cor:newton-poly}
If the slopes of $D$ are either all positive numbers or all negative numbers, the map
$1-\phi:D\to D$ is invertible.
\end{cor}

\begin{proof}
If the slopes of $D$ are all positive numbers, then $1+\phi+\phi^2+\cdots$
converges and is an inverse of $1-\phi$.

If the slopes of $D$ are all negative numbers, then the slopes of the dual isocrystal $D^{\vee}$
are all positive numbers. By choosing a basis of $D$, the matrix of $\phi^\vee$ is the transpose inverse of that of $\phi$, and $(1-\phi^{-t})=-\phi^{-t}(1-\phi^t)^{-1}=-\phi^{-t}(1+\phi^{t}+\phi^{2t}+\cdots)$ is invertible.
\end{proof}

\begin{defn}
Let $(F, \phi_F, \theta_F, N_F)$ be a filtered $(\phi, N)$-module with $P$-structure.
Consider the adjoint action $\bM \to \GL(\Lie U)$.
Assume $(F, \phi_F, \theta_F, N_F)\times^{P}\bM$
is weakly admissible.

The following conditions are equivalent:
\begin{enumerate}
\item 
The Hodge polygon of the $\phi$-module
$$(F, \phi_F, \theta_F, N_F)\times^{P}\bM\times^{\bM}\Lie U$$
has slopes $<0$;
\item
All Hodge-Tate weights of the semistable representation
$$
V_{\operatorname{ss}}((F, \phi_F, \theta_F, N_F)\times^{P}\Lie U)
$$ are negative integers,
where $V_{\operatorname{ss}}$ is the covariant Fontaine functor;
\item
$\theta_F^0(\Lie G) \cap \Lie U = \{0\}$.
\end{enumerate}
We say $\theta_F$ is $P^-$-dominant if
$(F, \phi_F, \theta_F, N_F)\times^{P}\bM$
is weakly admissible and
one of the conditions above holds.
\label{def:Pdom}
\end{defn}

\begin{lem}
\label{adj-phi-torsor}
If $(F, \phi_F, \theta_F, N_F)$ is weakly admissible
and $\theta_F$ is $P^-$-dominant,
then the map $1-\phi_{F\times^{P}\Lie U}$ 
is invertible.
\end{lem}

\begin{proof}
By Corallary \ref{cor:newton-poly}, it suffices to show all slopes of the Newton polygon of $F\times^P\Lie U$ are negative numbers,
which follows from the definition of
$P$-dominance and weak admissibility.
\end{proof}

\begin{lem}
\label{invertible-U}
If $(F, \phi_F, \theta_F, N_F)$ is weakly admissible
and $\theta_F$ is $P^-$-dominant,
then the map 
$$
U \to U,~u\mapsto u\phi_{F\times^{P}U}(u)^{-1}
$$
is invertible.
\end{lem}

\begin{proof}
Denote this map by $\Delta_F$.
Since $U$ is a unipotent algebraic group over
a ring of characteristic $0$,
the logarithm map
$\log: U\to \Lie U$ is an $\bM$-equivariant
bijection.
Let $m_0=m\in U$.
We inductively define
\begin{align*}
n_k & = (1-\phi_{F\times^{P}\Lie U})^{-1} \log m_k
\\
m_{k+1} & = n_k^{-1} m_k \phi_{F\times^PU}(n_k)
\end{align*}
Let $U_0=U \supset U_1 = [U,U] \supset \cdots\supset U_s=\{1\}$ be the lower central series.
It is clear that $m_k\in U_k$ for each $k=0,1,2,\cdots$.
Thus $m=\Delta_F(n_0n_1\cdots n_{s-1})$.
It remain to show the injectivity:
by the explicit formula, we conclude inductively on $k$ that
$x=n_k n_{k+1}\cdots n_{s-1}$
is the unique element of $U_k$
that solves the equation $\Delta_F(n_0\cdots n_{k-1})^{-1}m=\Delta_F(x)$.
\end{proof}

We didn't prove Lemma \ref{invertible-U}
by using the technique of Proposition \ref{prop:weak-adm-P}
because after embedding, we couldn't easily show the inverse actually lands
in $U$ (the argument would be equally lengthy).

\subsection{Framing}
Let $(F, \phi_F, \theta_F, N_F)$ be a filtered $(\phi, N)$-module.
Once we fix a (right) torsor morphism
$[-]:F\cong H(K_0\otimes \bC)$,
we have
$$
[\phi_F(h)]=[\phi_F]\varphi([h])
$$
for all $h\in H(K_0\otimes \bC)$
and some $[\phi_F]\in H(K_0\otimes \bC)$.
We remark that
changing framing as the effect of
changing $[\phi_F]$ to $ g^{-1}[\phi_F]\varphi(g)$,
$g\in H(K_0\otimes \bC)$.
The monodromy condition becomes
$N_F=p[\phi_F]\varphi(N_F)[\phi_F]^{-1}$
(recall that $N_F\in \Lie H$).

\begin{prop}
Fix a framing of $F\times^P\bM$.
If $(F, \phi_F, \theta_F, N_F)$ is weakly admissible
and $\theta_F$ is $P^-$-dominant,
then there exists a unique framing $[-]$ of $F$
such that $[\phi_F] = [\phi_{F\times^P\bM}]\in \bM(K_0\otimes\bC)$.
\label{prop:unique-phi-frame}
\end{prop}

\begin{proof}
Set $A:=[\phi_{F\times^P\bM}]$.
Suppose $[-]$ is a framing of $F$
that is compatible with the fixed framing of
$F\times^P\bM$, that is, 
$[\phi_F]=M A$, where $M\in U$.
Changing the framing $[-]$ amounts to change
$M A$ to $u^{-1}M A \varphi(u)$, $u\in U$.
We want to solve the equation
$$
u^{-1}M A \varphi(u) = A
\Leftrightarrow
M = u A \varphi(u)^{-1} A^{-1},
$$
which has a unique solution in $U$ by
Lemma \ref{invertible-U}.
\end{proof}

Fix a filtered $(\phi, N)$-module
$(\bar F, \phi_{\bar F}, \theta_{\bar F}, N_{\bar F})$
with $\bM$-structure.
As a consequence of 
Proposition \ref{prop:unique-phi-frame}
and \ref{prop:weak-adm-P},
the set of all possible monodromy operators
$N_F$ appearing in
a filtered $(\phi, N)$-module
$(F, \phi_F, \theta_F, N_F)$
such that 
$(F, \phi_F, \theta_F, N_F)^{\semis{\bM}}=
(\bar F, \phi_{\bar F}, \theta_{\bar F}, N_{\bar F})$
{\it does not depend on $(F, \phi_F, \theta_F, N_F)$}.
Indeed, the set of all possible monodromy operators
$N_{F}$ extending $N_{\bar F}$
naturally
form a $\bC$-vector space.

\begin{defn}
Denote by 
\[
\Ext(\WD(\bar F, \phi_{\bar F}, \theta_{\bar F}, N_{\bar F}), U):=
\left\{
N_F \ \middle| \
\begin{matrix}
\exists \text{ weakly admissible } (F, \phi_F, \theta_F, N_F)
\text{with $U\rtimes \bM$-structure}\\
\text{such that } (F, \phi_F, \theta_F, N_F)^{\semis{\bM}} = (\bar F, \phi_{\bar F}, \theta_{\bar F}, N_{\bar F})
\end{matrix}
\right\}.
\]
\end{defn}

\begin{lem}
We have
$$
\Ext(\WD(\bar F, \phi_{\bar F}, \theta_{\bar F}, N_{\bar F}), U)
=
\Ext(\WD(\bar F, \phi_{\bar F}, \theta_{\bar F}, N_{\bar F}), \Lie U).
$$
To elaborate,
abelian extensions to $\Lie U\rtimes \bM$
and non-abelian extensions to $U\rtimes \bM$
have the same (vector) space of monodromy operators.
\end{lem}

\begin{proof}
Clear by definition.
\end{proof}

We record the following standard fact for lack of  direct reference.

\begin{lem}
(1) If $P$ is a parabolic of a reductive group $H$,
then the center of $P$ is equal to the center of $H$, that is, $Z(P)=Z(H)$.

(2) If $M$ is a Levi subgroup of a reductive group $H$,
then the centralizer of the connected center of $M$ is $M$,
that is, $Z_H(Z(M)^\circ)=M$.
\end{lem}

\begin{proof}
(1) Say $P=U\rtimes M$, where $U$ is the unipotent radical.
First note $Z(P)\subset U\rtimes Z(M)$ because it is contained in $Z(M)$ after mod $U$.
Say $g = t u \in Z(P)$ where $t\in Z(M)$ and $U\in M$, and $s\in T$ where $T\subset M$ is a maximal torus.
We have $t u s u^{-1} t^{-1} = s\Leftrightarrow u=s u  s^{-1}$ for all $s\in T$, which forces $u=1$
by the root space decomposition of $U$.
So $Z(P)\subset Z(M)\subset T$. 
Let $t\in Z(P)$.
Say $U_{\alpha}\subset U$ is a root group,
where $\alpha\in X^*(T)$.
$t$ commuting with elements of $U_{\alpha}$
is equivalent to $t$ commuting with elements of $U_{-\alpha}$.
Therefore $Z(P)\subset Z_H(P^-)\cap Z_H(P)$
where $P^-$ is the opposite parabolic.
Finally since $P^-P$ is dense in $H$, we have
$Z(P)\subset Z(H)$.
Conversely, we have $Z(H)\subset T\subset P$
which implies $Z(H)\subset Z(P)$.

(2)
Indeed, $M=Z_H(\mu)$
for some cocharacter $\mu$.
\end{proof}

The following lemma also looks very standard
but we couldn't find a reference.

\begin{lem}
If $(F, \phi_F, \theta_F, N_F)$ is weakly admissible
and $\theta_F$ is $P^-$-dominant,
then $\theta_F\not=u \theta_F u^{-1}$
for all $u\in U$.
In particular, 
the element $u$ in Lemma \ref{triv-fil}
is unique.
\label{lemma:unique-u}
\end{lem}

\begin{proof}
Say $\theta_F = \lambda_{\can}$
for some $\lambda\in X_*(\bT)$
where $\bT\subset P$ is a maximal torus.
There exists another cocharacter $\mu\in X_*(\bT)$
such that $P=P_{\bG}(\mu)$.
Since both $P_{\bG}(\mu)$
and $P_{\bG}(\lambda)=\Aut^{\otimes}(\theta_F)$ are stable under the
$\bT$-action (by conjugation),
we can decompose
$\Lie (P \cap \Aut^{\otimes}(\theta_F))$
as a sum of $\bT$-weight spaces.
We claim that the unique
$\Gm$-weight of $\Lie (P \cap \Aut^{\otimes}(\theta_F))$ (where $\Gm$ acts
by $\mu:\Gm\to \bT$) is $0$:
indeed, $\Lie P$ has no negative weights
and $\Lie \Aut^{\otimes}(\theta_F)$ has no positive weights by the $P$-dominance.
Therefore $P \cap \Aut^{\otimes}(\theta_F) \subset Z_{\bG}(\mu)$ is contained in a maximal Levi subgroup of $P$.
So, $P \cap \Aut^{\otimes}(\theta_F) = P_{Z_{\bG}(\mu)}(\lambda)$.
By the previous lemma
(first use part 1, then use part 2),
we have $Z_{\bG}(\mu) = Z_{\bG}(Z(P \cap \Aut^{\otimes}(\theta_F))^\circ)$.

Now we prove this lemma.
Suppose $\theta_F=u \theta_F u^{-1}$
and $u\in U$.
Then 
$$
P \cap \Aut^{\otimes}(\theta_F) = P \cap \Aut^{\otimes}(u \theta_F u^{-1})
= u P \cap \Aut^{\otimes}(\theta_F) u^{-1}
$$
which implies
$Z_{\bG}(\mu)=uZ_{\bG}(\mu)u^{-1}$
$\Rightarrow$
$u\in Z_{\bG}(\mu)$
$\Rightarrow$
$u=1$.
\end{proof}

\begin{cor}
Let $(\bar F, \phi_{\bar F}, \theta_{\bar F}, N_{\bar F})$
be a weakly admissible filtered $(\phi, N)$-module
with $\bM$-structure with $\theta_{\bar F}$ being $P^-$-dominant.
The set of isomorphism classes
of weakly admissible filtered $(\phi, N)$-modules $(F, \phi_{F}, \theta_{F}, N_F)$ equipped with
an identification
$(F, \phi_{F}, \theta_{F}, N_F)\times^P\bM 
=  (\bar F, \phi_{\bar F}, \theta_{\bar F}, N_{\bar F})$
is in canonical bijection with
the cartesian product
$$
U(K\otimes \bC) \times 
\Ext(\WD(\bar F, \phi_{\bar F}, \theta_{\bar F}, N_{\bar F}), U).
$$
\label{cor:U-torsor-extn}
\end{cor}

\begin{proof}
It is harmless to fix a framing of $(\bar F, \phi_{\bar F}, \theta_{\bar F}, N_{\bar F})$.
By Proposition 
\ref{prop:unique-phi-frame},
each isomorphism class
of such filtered $(\phi, N)$-module with $P$-structure
admits a unique framing
under which
$[\phi_F]=[\phi_{\bar F}]$.
By Proposition \ref{prop:weak-adm-P},
the choice of $\theta_F$
and $N_F$ are independent.
By
Lemma \ref{triv-fil}
and Lemma 
\ref{lemma:unique-u},
the choice of $\theta_F$
is parameterized by $U(K\otimes C)$.
\end{proof}

\begin{cor}
Let $U_1\unlhd U_2$ be closed algebraic subgroups
of $U$ which are stable under adjoint action of $\bM$.
Corollary \ref{cor:U-torsor-extn} holds
after replacing $U$ by $U_1/U_2$.
\label{cor:U-torsor-extn+}
\end{cor}

\begin{proof}
Clear.
\end{proof}

\begin{thm}
\label{thm:semistable}
Let $\rho:\Gal_K\to \bM(\bC)$
be a semistable Galois representation
with $P^-$-dominant Hodge type.
Let $\Gal_K$
act on $\Lie U$ via
$\Gal_K\xrightarrow{
\rho
}\bM(\bC)\xrightarrow{\Ad} \Aut(U)$.
Let $U_1\unlhd U_2$ be closed algebraic subgroups
of $U$ which is stable under adjoint action of $\bM$.
Set $U':=U_2/U_1$.

(1)
There is a canonical bijection
$$
H^1(\Gal_K, U'(\bC))
\cong H^1(\Gal_K, \Lie U'(\bC))
$$
which is compatible with the long exact sequence
induced from the short exact sequence
$U_1\to U_2\to U'$.

(2)
There is a commutative diagram
$$
\xymatrix{
H^1(\Gal_K, U(\bC)) \ar[r]^{\cong} \ar[d] & H^1(\Gal_K, \Lie U(\bC)) \ar[dl] \\
\Ext(\WD(\bar F, \phi_{\bar F}, \theta_{\bar F}, N_{\bar F}), U) \
}
$$
with all fibers canonically identified with
$U(K\otimes \bC)$.

(3)
There is a isomorphism
of $\bC$-vector spaces
$$
\Ext(\WD(\bar F, \phi_{\bar F}, \theta_{\bar F}, N_{\bar F}), U)
\cong
H^2(\Gal_{K}, \Lie U(\bC)).
$$
\end{thm}

\begin{proof}
Part (1, 2): by Corollary \ref{cor:U-torsor-extn+},
we specify a bijection
$$
U(K\otimes \bC) \times 
\Ext(\WD(\bar F, \phi_{\bar F}, \theta_{\bar F}, N_{\bar F}), U)
\to
\Lie U(K\otimes \bC) \times 
\Ext(\WD(\bar F, \phi_{\bar F}, \theta_{\bar F}, N_{\bar F}), U).
$$
Set $(u, N)\mapsto (\log u, N)$
where $\log$ is automatically
truncated after $\operatorname{nc}(U)$ terms
(the nilpotency class of $U$).
This shows
$
H^1_\st(\Gal_K, U'(\bC))
\cong H^1_\st(\Gal_K, \Lie U'(\bC))
$
where $H^1_\st\subset H^1$
is the subset of semistable extensions.
By \cite[Lemma 6.4, Lemma 6.5]{Be02},
$
H^1_\st(\Gal_K, \Lie U'(\bC))
\cong H^1(\Gal_K, \Lie U'(\bC))
$,
and so we have a commutative diagram
\begin{align*}
\raisebox{-4ex}{\text{($*$)}}
&
\xymatrix{
H^1_{\mathrm{st}}(\mathrm{Gal}_K, U'(\mathbb{C})) \ar[r]^{\cong} \ar@{^{(}->}[d] & H^1_{\mathrm{st}}(\mathrm{Gal}_K, \mathrm{Lie}\, U'(\mathbb{C})) \ar[d]^{\cong} \\
H^1(\mathrm{Gal}_K, U'(\mathbb{C}))  & H^1(\mathrm{Gal}_K, \mathrm{Lie}\, U'(\mathbb{C}))
}
\end{align*}
It remains to show the left vertical map is surjective,
and we do this by considering the long exact sequence
$$
\xymatrix{
  H^1_{\mathrm{st}}(\mathrm{Gal}_K, U_1(\mathbb{C})) \ar[d] &
  H^1_{\mathrm{st}}(\mathrm{Gal}_K, U_2(\mathbb{C})) \ar[d]&
  H^1_{\mathrm{st}}(\mathrm{Gal}_K, U'(\mathbb{C})) \ar[d] \\
H^1(\mathrm{Gal}_K, U_1(\mathbb{C})) \ar[r] & H^1(\mathrm{Gal}_K, U_2(\mathbb{C})) \ar[r] & H^1(\mathrm{Gal}_K, U'(\mathbb{C}))
}
$$
we argue by induction on $\dim U'$
and so we assume the left and right vertical arrows are surjective; it is clear the middle vertical arrow
is also surjective
using diagram ($*$).
We warn the reader that a priori there are no horizontal arrows in the top arrow.

Part (3):
since $\Lie U(\bC)$
has negative Hodge-Tate weights,
we have $H^0(\Gal_K, \Lie U(\bC))=0$.
By the Euler characteristic,
$\dim H^1 = \dim H^0+\dim H^2+\dim U(\bC)
=\dim H^2+\dim U(\bC)$.
So, part (3) follows from part (2)
by noting the right vertical map is linear.
\end{proof}

\begin{cor}
\label{cor:semistable}
Let $\rho: \Gal_{K}\to \bM(\bC)$ be a semistable representation whose Hodge type is $P^-$-dominant.
Endow $\Lie U(\bC)$ with $\Gal_K$-action 
via $\ad\rho$.

Any $\rho':\Gal_K\to P(\bC)$
such that $\rho'\times^P\bM= \rho$
is semistable
with its possible Weil-Deligne
types parameterized by $H^2(\Gal_K, U(\bC))$.
\end{cor}

\begin{proof}
It is already proved in the proof of Theorem \ref{thm:semistable}.
Alternatively, we can deduce it from
Theorem \ref{thm:semistable}
combined with Corollary \ref{cor:U-torsor-extn+}.
\end{proof}

\begin{cor}
\label{cor:crys}
Let $\rho: \Gal_{K}\to \bM(\bC)$ be a crystalline
representation whose Hodge type is $P^-$-dominant.
Endow $\Lie U(\bC)$ with $\Gal_K$-action 
via $\ad\rho$.

Assume 
$H^2(\Gal_K, U(\bC))=0$.
Then any $\rho':\Gal_K\to P(\bC)$
such that $\rho'\times^P\bM= \rho$
is crystalline.
\end{cor}

\begin{proof}
By Corollary \ref{cor:semistable},
the Weil-Deligne type of $\rho'$
is unique.
But $\rho$ itself 
composed with $\bM\hookrightarrow P$
is a crystalline representation.
\end{proof}

\section{Weil-Deligne types of de Rham extensions}

\subsection{Non-abelian Shapiro lemma}
Let $G$ be a possibly non-split reductive group over $K$
that splits over $L$
and let $\lsup LG:=\bG\rtimes \Gal(L/K)$.
There is an equivalence of groupoids
$$
\{\text{parameters
$\Gal_K\to \lsup LG$}\}
\leftrightarrow
\{\text{parameters
$\Gal_{\bQ_p}\to \lsup L\Res_{K/\bQ_p}G$}\},
$$
see, for example, \cite[Lemma 9.4.1]{GHS} for unramified
groups and \cite{St10} for the general case.
As a consequence,
we will assume $K=\bQ_p$ from now on
as all results for $\bQ_p$
automatically imply
the same statement for general $K$.

The phrases ``Galois representation''
and ``parameters''
are used interchangeably.

Set $\lsup LP:=P\rtimes \Gal(L/\bQ_p)$
where $P=U\rtimes \bM$ is a $\Gal(L/\bQ_p)$-stable
parabolic of $\bG$.

\subsection{Weil-Deligne types}
\label{subsec:WD}
Let $\rho:\Gal_{\bQ_p}\to \lsup LP(\bC)$ be a de Rham
Galois representation.
We choose
the splitting field $L/\bQ_p$ of $G$
large enough so that
$\rho$ becomes semistable
after restriction to $\Gal_L$.
We denote by $\rho|_{\Gal_L}:\Gal_L\to \bG(\bC)$
the restriction of $\rho$ composed with
the evaluation at identity $\lsup LG\to \bG$.

$D_{\sem}^L(\rho|_{\Gal_L})=:
(\phi_{\rho}, \theta_{\rho}, N_{\rho})$
is a weakly admissible filtered $(\phi, N)$-module
that admits a semilinear action
of $\Gal(L/\bQ_p)$
which we denote by $\tau_{\rho}$.
We forget the filtration structure
and define
$\WD(\rho):=
(\phi_{\rho}, N_{\rho}, \tau_\rho)$
and call it {\it the Weil-Deligne type}
of $\rho$.
Our definition is consistent with 
\cite[Definition 2.6.1]{BG19}.

\subsection{Remark}
It is possible to remove the dependence on the splitting field
$L$ by descending $\WD(\rho)$
to the $\sigma$-isotypic part
$\WD(\rho)_{\sigma}$
for a fixed embedding $\sigma:L_0\hookrightarrow\bC$
and combine the data
of $\phi_{\rho}$
and $\tau_{\rho}|_{I_{\bQ_p}}$
into a representation
of the Weil group $W_{\bQ_p}:=\bZ\ltimes I_{\bQ_p}$
to get a genuine Weil-Deligne representation
$(r_\sigma, N_\sigma)$
on 
$\WD(\rho)_{\sigma}$
(c.f. \cite[Lemma 2.6.6]{BG19}).
Moreover, if $\rho$ is potentially crystalline
and $G$ is an unramified group,
then 
$(r_\sigma, N_\sigma)$
can be reconstructed from
$r_\sigma|_{I_{\bQ_p}}$ alone,
which is often called {\it the inertial type}
of $\rho$ in the literature.
In the rest of this note,
we reserve the notation $\WD_\sigma(\rho)$
for the inertial type $r_\sigma|_{I_{\bQ_p}}$
and we remark that
\cite{EG23}
uses the notation $\WD(\rho)$
for inertial types,
not Weil-Deligne types.

\begin{thm}
\label{thm:de-rham}
Let $\rho:\Gal_{\bQ_p}\to \lsup LM(\bC)$ be a de Rham
Galois representation
of $P^-$-dominant Hodge type.
Let $\Gal_{\bQ_p}$ act on $U(\bC)$
by $\Gal_{\bQ_p}\xrightarrow{\rho}\lsup LM(\bC)\xrightarrow{\Ad}U(\bC)$.
Then
there is a commutative diagram
$$
\xymatrix{
H^1(\Gal_{\bQ_p}, U(\mathbb{C})) \ar[r]^\sim 
\ar@{=}[d] & H^1(\Gal_{\bQ_p}, \operatorname{Lie} U(\mathbb{C})) \ar@{=}[d] \\
H^1(\Gal_L, U(\mathbb{C}))^{\Gal(L/\bQ_p)} \ar[r]^\sim & H^1(\Gal_L, \operatorname{Lie} U(\mathbb{C}))^{\Gal(L/\bQ_p)}.
}
$$
Here $L$ is chosen so that
$\rho|_{\Gal_L}$ is semistable.
\end{thm}

\begin{proof}
The vertical maps are equality by
\cite[Theorem 3.15]{Ko02}.
It remains to show
that the isomorphisms in
Theorem \ref{thm:semistable}
are compatible with $\Gal(L/\bQ_p)$-action,
which is clear.
\end{proof}

\begin{cor}
Theorem \ref{thm:semistable}
and Corollary \ref{cor:semistable}
both hold after replacing ``semistable''
by ``de Rham'',
$P$ by $\lsup LP$
and $\bM$ by $\lsup LM$.
\label{cor:de-rham}
\end{cor}

\begin{proof}
Clear.
\end{proof}

\section{Moduli of de Rham parameters: the parabolic case}

By \cite{Min24},
the Emerton-Gee stacks
$\cX_{K, \lsup LP}$ (or $\cX_{\lsup LP}$
when the choice of $K$ is clear from the context)
is a formal algebraic stack, locally
Noetherian
over $\spf \bar\bZ_p$.

\subsection{Remark}
When $\lsup LP$ is connected,
then there is no ambiguity in $\cX_{K, \lsup LP}$.
When $\lsup LP$ is disconnected,
we always use the version such that
$\cX_{K, \lsup LP}(\bFp)$
is the groupoid
of parameters $\Gal_K\to \lsup LP(\bFp)$,
up to $P$-conjugacy rather than $\lsup LP$-conjugacy.

\begin{lem}
$\cX_{K, \lsup LP}$
is a Noetherian formal algebraic stack over $\spf \bar\bZ_p$.
\end{lem}

\begin{proof}
It suffices to show $\cX_{K, \lsup LP, \red}$
is quasi-compact,
which follows from 
\cite[Proposition 10.1.8]{L23B}.
\end{proof}

\subsection{$\bar\bZ_p$-flat quotient of Noetherian formal
algebraic stacks}
Let $\cY$ be a Noetherian formal algebraic stack
over $\spf \bar\bZ_p$.
Let $\spf A$ be a Noetherian smooth cover of $\cY$,
which exists by \cite[Tag 0AIC, 0AID]{Stacks}
and the definition of formal algebraic stacks
\cite[Definition 5.3]{Eme}.
Let $T$ be the $p$-power torsion ideal of $A$
and define the $\bar\bZ_p$-flat quotient of $\spf A$
to be $\spf A/T$.
Such a construction descends along smooth covers.

\begin{defn}
Let $\lambda$ be a $P^-$-dominant Hodge type
and let $\tau$ be an inertial type,
both valued in $\lsup LP$.

Define
$\cX_{K, \lsup LP}^{\tau, \lambda}$
to be (if it exists)
the Noetherian formal algebraic stack
which is flat over $\spf \bar\bZ_p$,
uniquely determined as
the $\cO$-flat closed substack of $\cX_{K, \lsup LP}$
by the following property:
if $A^\circ$ is a finite $\bar\bZ_p$-flat algebra,
then  
$\cX_{K, \lsup LP}^{\tau, \lambda}(A^\circ)$
is the subgroupoid consisting of $L$-parameters
which become potentially semistable
of Hodge type $\lambda$
and inertia type $\tau$ after inverting $p$.
\label{def:XP}
\end{defn}

We record the following obvious properties:
\begin{lem}
(1)
If
$\cX_{K, \lsup LP}^{\tau, \lambda}$
exists,
then it is indeed unique.

(2)
If $G=G_1\times G_2$ is reductive,
then
$$
\cX_{K, \lsup L(G_1\times G_2)}^{(\tau_1, \tau_2),
(\lambda_1, \lambda_2)}
=
\cX_{K, \lsup LG_1}^{\tau_1,
\lambda_1}
\times
\cX_{K, \lsup LG_2}^{\tau_2,
\lambda_2}
$$

(3)
If $i:G\to H$ is a central isogeny
of reductive groups (with finite kernel),
then
$\cX_{K, \lsup LH}^{\tau, \lambda}\to\cX_{K, \lsup LG}^{i(\tau),i(\lambda)}$
is relatively representable by a finite morphism.

(4, Shapiro Lemma)
Under the identification (\cite[Proposition 7.2.4]{L23B})
$$
\Sha: \cX_{K, \lsup LG} \cong \cX_{\bQ_p, \lsup L{\Res_{K/\Qp}G}},
$$
we have
$$
\cX_{K, \lsup LG}^{\tau, \lambda}
=
\cX_{\bQ_p, \lsup L{\Res_{K/\Qp}G}}^{\Sha(\tau),~\Sha(\lambda)}.
$$

(5, the tame case)
If $G$ is tame over $K$ and reductive,
then $\cX_{K, \lsup LG}^{\tau, \lambda}$
exists
and is $p$-adic formal algebraic
of finite type over $\spf \bar\bZ_p$.
\label{lem:XP}
\end{lem}

\begin{proof}
(1)
Suppose $\cX_1$ and $\cX_2$
are closed substacks satisfy the defining properties.
So does $\cX_3:=\cX_1\times_{\cX_{K, \lsup LP}}\cX_2$.
Both $\cX_3\to \cX_1$ and $\cX_3\to \cX_2$
are formally \'etale closed immersions.
So $\cX_3$ is an open and closed substack of $\cX_1$;
if $\cX_1$ has a component disjoint from $\cX_3$,
then $\cX_3$ cannot satisfy the defining property.

(2) and (4) Both follows immediately from the uniqueness.

(3)
Say $X_1$ and $X_2$
are two objects of $\cX_{K, \lsup LH}^{\tau, \lambda}(\spec A)$
that map to the same object of
$\cX_{K, \lsup LG}^{i(\tau), i(\lambda)}(\spec A)$,
then it makes sense to take
their ratio ``$X_1/X_2$'',
which is a well-defined object
of $\cX_{K, Z}$
where $Z$ is the kernel of $\bH\to \bG$.
Since $Z$ is a multiplicative finite group scheme,
$\cX_{K, Z}$
is the same as the moduli
of Weil representations
by  \cite[Corollary 2.7.4]{L23B}
and is finite over $\spf \bar\bZ_p$.

(5)
It is proved in the Appendix \ref{sec:moduli-ss}.
See Theorem \ref{thm:ss}.
\end{proof}

\begin{thm}
If $p>3$,
then 
$\cX_{K, \lsup LP}^{\tau, \lambda}$
exists,
and
$\cX_{K, \lsup LG}^{\tau, \lambda}$
is of finite type over $\spf\bar\bZ_p$.
\label{thm:XP}
\end{thm}

\begin{proof}
We first assume $\lsup LP=\lsup LG$ is reductive.
Let $G'$ be the simply-connected cover
of $G$ and let $Z$ be the kernel of $G$.
Then $G'\times Z\to G$ is a central isogeny.
Since $p>3$, $G'$ is a finite product
of Weil restrictions of tame groups.
So, $\cX_{K, \lsup LP}^{\tau, \lambda}$ exists
by Lemma \ref{lem:XP}.

In general, write $i$ for the surjective homomorphism
$\lsup LP\to\lsup LM$
and set
$
\cX:=\cX_{K, \lsup LP} \underset{\cX_{K, \lsup LM}}{\times}
\cX_{K, \lsup LM}^{i(\tau), i(\lambda)}.
$
It is clear $\cX_{K, \lsup LP}^{\tau, \lambda}$
is representable by the $\bar\bZ_p$-flat quotient of $\cX$,
as a consequence of Corollary \ref{cor:de-rham}.
\end{proof}

\begin{defn}
Fix a Hodge type $\lambda$.
Define $\cX_{K, \lsup LP}^\lambda$
to be the scheme-theoretic union
of the finitely many
$\cX_{K, \lsup LP}^{\lambda, \tau}$
as $\tau$ ranges over all inertial types.
$\cX_{K, \lsup LP}^\lambda$
is a $p$-adic formal algebraic stack.
\end{defn}

\section{Borel-valued Weil-Deligne moduli stacks}

Without loss of generality, set $K=\bQ_p$.

\subsection{The moduli of Weil-Deligne representations}
$\WD_{\lsup LB}$
is constructed in \cite[Definition 2.1.2]{BG19}
and is denoted by $[\cY_{L/K, \varphi, \cN}/\bG]$ 
following the notation of loc. cit..
We reiterate that
for any fixed choice of Hodge type $\lambda$,
the usual linear version of Weil-Deligne representations
is equivalent to the semilinear version
of Weil-Deligne representations
considered in Subsection \ref{subsec:WD}
by the explicit formula in \cite[Lemma 2.6.6]{BG19}.

\begin{prop}
\label{prop:sp}
(1)
$\cX^{\rho}_{\lsup LB}[\frac{1}{p}]\to
\WD_{\lsup LB}$
is surjective.

(2)
Part (1) of Theorem \ref{thm:sp} holds,
and $\WD_{\lsup LB}$
is equidimensional of dimension $0$
over $\bar\bQ_p$.
\end{prop}

\begin{proof}
(1)
If $(\phi_F, N_F, \tau_F)$
is a Weil-Deligne type valued
in $\lsup LB$,
then $(\phi_F, \rho, N_F)$
is a weakly admissible filtered $(\phi, N)$-module
with Galois descent data furnished by $\tau_F$.

(2) By \cite[Lemma 4, Lemma 14]{L25} and part (1),
it suffices to show
all fibers of 
$\cX^{\rho}_{\lsup LB}[\frac{1}{p}]\to
\WD_{\lsup LB}$
are irreducible of constant rank.
Indeed for any $\spec \bQp\to \WD_{\lsup LB}$,
the fiber is kernel
to a surjective linear map
$H^1(\Gal_K, \Lie U(\bQp))\to H^2(\Gal_K, \Lie U(\bQp))$
by Theorem \ref{thm:de-rham}
and Corollary \ref{cor:de-rham},
which is of constank rank $[K:\bQ_p]\dim U$.
\end{proof}

\subsection{Remark}
The morphism $\cX^{\rho}_{\lsup LB}[\frac{1}{p}]\to
\WD_{\lsup LB}$
is literally the forgetful functor
that forgets the filtration structure
on a weakly admissible
filtered $(\phi,N)$-module
with descent data and $\bG$-structure.

\subsection{}
Next, we analyze the structure of $\WD_{\lsup LB}$.
Inspired the discussion in \cite[Section 10]{L25},
we decompose $\WD_{\lsup LB}$
according to the set of simple positive roots
$\Delta=\{\alpha_1, \dots, \alpha_s\}
=\Delta_G(B)$.

\begin{lem}
(1)
$
\cX^{\rho}_{\lsup LB}[\frac{1}{p}]
\to
\cX^{\rho}_{\prod_{i=1}^sU_{\alpha_i}\rtimes\lsup LT}[\frac{1}{p}]
$
is surjective with
all fibers representable by an affine space
of dimension $(\dim U - s)$.

(2)
We have
$$
\WD_{\lsup LB} \cong \WD_{\prod_{i=1}^sU_{\alpha_i}\rtimes \lsup LT}.
$$
\label{lem:comb}
\end{lem}

\begin{proof}
Note that 
$H^2(\Gal_{\bQ_p}, \Lie U(\bQp))\cong H^2(\Gal_{\bQ_p}, \prod_{i=1}^sU_{\alpha_i}(\bQp))$
when the Hodge type is $\rho$.

(1)
It follows from Theorem \ref{thm:semistable}
and Theorem \ref{thm:de-rham}.

(2)
Repeating the proof of Proposition
\ref{prop:sp},
we have
$\cX^{\rho}_{\prod_{i=1}^sU_{\alpha_i}\rtimes\lsup LT}[\frac{1}{p}]\to
\WD_{\prod_{i=1}^sU_{\alpha_i}\rtimes\lsup LT}$
is also surjective with all fibers
being representable by the kernel
of
$H^1(\Gal_{\bQ_p}, \prod_{i=1}^sU_{\alpha_i}\rtimes\Lie T(\bQp))\to H^2(\Gal_{\bQ_p}, \Lie U(\bQp))$.

The upshot is 
$\WD_{\lsup LB} \to \WD_{\prod_{i=1}^sU_{\alpha_i}\rtimes \lsup LT}$
is surjective
with all fibers representable by a point,
and is thus an isomorphism.
\end{proof}

Let $L_1/\bQ_p$ be the splitting field
of $\alpha_1$,
and set
$$
\Delta^1:=\Delta \backslash \{\gamma\cdot\alpha_1|
\gamma\in \Gal(L_1/\bQ_p)\}.
$$
Write $U_{\alpha_i}\subset U$ for the
root group of weight $\alpha_i^\vee$.
Let $\bM_1\subset \bG$ be the
Levi subgroup generated by the root groups
$U_{\pm\alpha}$, $\alpha\in\Delta^1$.
Set $\lsup LM_1:=\bM_1\rtimes \Gamma$.
Then $\lsup LB_1:=\lsup LB\cap \lsup LM_1$
is a Borel of $\lsup LM_1$.
We have
\begin{align*}
\WD_{\lsup LB} & \cong \WD_{\prod_{\alpha\in \Delta}\rtimes \lsup LT}\\
\WD_{\lsup LB_1} & \cong \WD_{\prod_{\alpha\in \Delta^1}\rtimes \lsup LT}
\end{align*}
Since 
$
\prod_{\alpha\in \Delta^1}\rtimes \lsup LT
\hookrightarrow
\prod_{\alpha\in \Delta}\rtimes \lsup LT
\to
\prod_{\alpha\in \Delta^1}\rtimes \lsup LT
$
composes to identity,
we conclude that 
$\Ch_0(\WD_{\lsup LB_1})$
is a direct summand of 
$\Ch_0(\WD_{\lsup LB})$.
Set
$\WD_{\lsup LB}^{\alpha_1-st}$
to be the closure of the complement
of 
$\WD_{\lsup LB_1}$
in
$\WD_{\lsup LB}$
where ``st'' means semi-stable.

\begin{lem}
All $\bQp$-fibers of 
$\WD_{\lsup LB}^{\alpha_1-st}
\to \WD_{\lsup LB_1}$
are isomorphic to 
$$H^2(\Gal_{\bQ_p}, \prod_{\gamma\in \Gal(L_1/\bQ_p)}U_{\gamma\cdot\alpha_1}(\bQp))\cong H^2(\Gal_{L_1}, U_{\alpha_1}(\bQp)).$$
\end{lem}

\begin{proof}
It follows from the Shapiro's Lemma.
\end{proof}

In particular, irreducible components of
$\WD_{\lsup LB}^{\alpha_1-st}$
are in natural bijection 
that of its image in 
$\WD_{\lsup LB_1}$,
by \cite[Lemma 14]{L25}.
We will be satisfied with dealing with
a special case:

\begin{prop}
\label{prop:main}
If $M_1$ has a direct summand $S\cong\Res_{L_1/\bQ_p}\Gm$
and the image of $\gamma\cdot\alpha_1$ ($\gamma\in \Gal(L_1/\bQ_p)$)
in $X^*(S)$
form a basis of
$X^*(S)$,
then
$\WD_{\lsup LB}^{\alpha_1-st}
\to \WD_{\lsup LB_1}
\to \WD_{\lsup L(M_1/S)}$
induces an isomorphism
of top Chow groups.

In particular, Theorem \ref{thm:main} holds.
\end{prop}

\begin{proof}
The composite $\WD_{\lsup LB}^{\alpha_1-st}
\to \WD_{\lsup L(M_1/S)}$
is surjective with all fibers being
$H^2(\Gal_{L_1}, U_{\alpha_1}(\bQp))$.
So,
$\Ch_0(\WD_{\lsup LB})\cong
\Ch_0(\WD_{\lsup LB_1})\oplus \Ch_0(\WD_{\lsup L(M_1/S)})$.

The ``in particular'' part follows
immediately from \cite[Proposition 7]{L25}.
\end{proof}

\subsection{The potentially crystalline part}

Denote by $\WD_{\lsup LB}^{\pcrys}\subset \WD_{\lsup LB}$
the substack where the monodromy operator vanishes.
The projection
$\WD_{\lsup LB}^{\pcrys}\to \WD_{\lsup LT}$
is clearly an isomorphism.
By the local Langlands duality for algebraic tori,
langlands parameters (Weil form)
$W_{\bQ_p}\to\lsup LT(\bQp)$
are equivalent to 
smooth
characters
$T(\bQ_p)\to\bQp^{\times}$.

We recall the structure theory of $p$-adic tori
in general.
Let $T$ be an algebraic tori over $\bQ_p$.
$T(\bQ_p)$ has a maximal bounded subgroup,
denoted by $T(\bQ_p)^1$
and an Iwahori subgroup, denoted by $T(\bQ_p)^0$
(\cite[Section 2.5(c)]{KP}).
The quotient
\begin{equation}
\label{eq:kot}
T(\bQ_p)^1/T(\bQ_p)^0 \cong X_*(T)_{I_{\bQ_p}, \tor}^{\Gamma}
\end{equation}
and
$$
T(\bQ_p)/T(\bQ_p)^0 \cong X_*(T)_{I_{\bQ_p}}^{\Gamma},
$$
c.f. \cite[Section 11.7]{KP}.

\begin{lem}
\label{lem:kottwitz}
Assume $p>2$.
The coarse moduli space of $\WD_{\lsup LT}$
is isomorphic to the moduli stack of
smooth characters
$\spec A\mapsto \Hom(T(\bQ_p), A^\times)$,
and part (4) of Theorem \ref{thm:sp} holds.
\label{thm:coarse}
\end{lem}

\begin{proof}
Denote by $K$ the maximal prop-$p$ subgroup
of $T(\bQ_p)^0$.
Any smooth homomorphism
$K\to \bQp^\times$
has finite pro-$p$ image
but the torsion part of $\bQp^\times$
is prime-to-$p$.
So any smooth character factors through
$$
T(\bQ_p)/K \cong 
X^*(T)_{I_{\bQ}, \tf}^\Gamma
\times T(\bQ_p)^1/K.
$$
by the structure of finitely generated abelian groups.
Since $X^*(T)_{I_{\bQ}, \tf}^\Gamma\cong \bZ^{\oplus r}$
is a finiely generated free abelian group,
we have
$|\WD_{\lsup LT}|\cong \bA^{\oplus r}
\times \{\text{characters of the torsion group
$T(\bQ_p)^1/K$}\}$.
\{\text{characters of the torsion group
$T(\bQ_p)^1/K$}\}
is representable by the
constant group scheme of the Pontryagin dual
of $T(\bQ_p)^1/K$,
and we have a short exact sequence
$$
1\to T(\bQ_p)^0/K \to T(\bQ_p)^1/K
\to T(\bQ_p)^1/T(\bQ_p)^0\to 1.
\qedhere
$$
\end{proof}

\section{The cyclotomic-free case}

We still assume $K=\bQ_p$.
Denote by $\cX:=\cX_{\lsup LG, \red}$
the reduced Emerton-Gee stack.
Since $\rho\mapsto H^2(\Gal_{\bQ_p}, \ad(\rho))$ is a coherent sheaf
over $\cX$,
the vanishing locus of $H^2(\Gal_{\bQ_p}, \ad(-))$
is an open substack of $\cX$,
which we denote by $\cX^{\cf}=\cX^{\cf}_{\lsup LG}$.

In order to be able to use results of \cite{L23A},
we assume either $G$ is tamely ramified.

\begin{lem}
Each $\bar\rho\in \cX(\bFp)$
factors through $U\rtimes\lsup LS$,
where $S$ is a $\bQ_p$-rational maximal torus
of $G$, $\lsup LS$ is the Langlands dual group of $S$,
and $U$ is the unipotent radical of a parabolic
$\lsup LP\supset \lsup LS$
subgroup of $\lsup LG$.
Moreover, we can demand the splitting field of $S$
to be unramified over
$L$.
\label{lem:ls}
\end{lem}

The tori $S$ in Lemma \ref{lem:ls}
are usually called {\it maximally unramified maximal tori}
of $G$.

\begin{proof}
It is the main theorem of \cite{L23A}.
\end{proof}

\begin{lem}
\label{lem:dim-Xcf}
(1) 
There exists a finite set of types
$\Xi=\{(\lambda, \tau)\}$, such that
each $\bar\rho \in \cX^{\cf}(\bFp)$
admits a potentially crystalline lift of
regular Hodge type $\lambda$
and inertial type $\tau$
with $(\lambda, \tau)\in \Xi$.

(2) $\cX^{\cf}$ is equidimensional of dimension
$\dim \bG/\brB$.
\end{lem}

\begin{proof}
We first explain the construction of $\Xi$.
For each stable conjugacy class
of maximally unramified maximal tori $S$ of $G$
that splits over $L_S$,
choose an arbitrary regular cocharacter 
$\lambda_S\in X_*(\bS)\subset X_*(\bG)$.
The set $\Xi$ consists of pairs $(\lambda, \tau)$ 
where $\lambda=\lambda_S$ for some $S$
and $\tau$ is an arbitrary inertial type
that splits over $L_S$.

(1)
By Lemma \ref{lem:ls}, 
$\bar\rho^{\sem}$ factors through
$\lsup LS(\bFp)$ for some $S$,
which by the LLC corresponds to
a character $\chi:S(\bQ_p)\to \bFp^\times$.
By Lemma \ref{lem:exist-de-rham}, there exists
a potentially crystalline character
$\chi_\lambda:S(\bQ_p)\to \cO_{\bC}^\times$
of Hodge type 
$\lambda=\lambda_S$.
Take the Teichmuller lift
$[\chi \chi_\lambda^{-1}]:S(\bQ_p)\to \cO_{\bC}^\times$
and multiply it by $\chi_\lambda$,
we get a lift of $\chi$ of Hodge type $\lambda$
and inertial type that splits over $L_S$.

Say  $\bar\rho$
factors through a parabolic $\lsup LP =:U \rtimes \lsup LM$.
Note that $\bar\rho$ as an extension of $\bar\rho^{\sem}$
is represented by a class $[\bar c]$ in $H^1(\Gal_{\bQ_p}, U(\bFp))$.
Since $H^2(\Gal_{\bQ_p}, U(\bFp))=0$,
 $H^1(\Gal_{\bQ_p}, U(\bZp))\to H^1(\Gal_{\bQ_p}, U(\bFp))$
 is surjective,
 and therefore there exists a class $[c]$ lifting $[\bar c]$.
 By Theorem \ref{thm:de-rham},
 $[c]$ represents a de Rham lift of $\bar\rho$
 that satisfies the desired properties.

 (2)
 The finiteness of $\Xi$ implies
$\mathcal{X}$ is equal to
 the reduction of the union
 of finitely many potentially crystalline stacks,
 of regular Hodge type.
 The lemma follows from the equidimensionality
 of potentially crystalline stacks,
 of regular Hodge type.
 See Appendix \ref{sec:moduli-ss}.
\end{proof}

\begin{lem}
Let $U$ be a unipotent algebraic group
which admits an action of $\lsup LG$.
Then $\cX^{\cf}_{U\rtimes \lsup LG}$
has dimension $\dim U + \dim\bG/\brB$.
\label{lem:dim-UG}
\end{lem}

\begin{proof}
It is a trivial extension of Lemma \ref{lem:dim-Xcf}.
\end{proof}

\begin{cor}
The scheme-theoretic image of
$\cX_{\lsup LB}^{\cf}$
in $\cX^{\cf}$
is the entire $\cX^{\cf}$.
\end{cor}

\begin{proof}
By Lemma \ref{lem:ls},
$\cX^{\cf}$ is the scheme-theoretic union
of $\cX^{\cf}_{U\rtimes\lsup LS}$
where $S$ ranges through
the (finitely many stable conjugacy classes)
of $\bQ_p$-rational maximal tori of $G$.
By Lemma \ref{lem:dim-UG},
$\dim \cX^{\cf}_{U\rtimes\lsup LS} < \dim \cX^{\cf}$
unless $U$ is the unipotent radical of $\lsup LB$
and therefore $\lsup LS\subset \lsup LB$,
up to $\bG$-conjugacy.
\end{proof}

In the rest of this section,
$U$ is specialized to be
the unipotent radical of $\lsup LB$,
and $\lsup LS$
is specialized to $\lsup LT$.

\begin{lem}
The morphism
$\cX^{\cf}_{\lsup LB}\to \cX^{\cf}_{\lsup LT}$
induces a bijection of irreducible components.
\label{lem:LB-LT}
\end{lem}

\begin{proof}
We replace $\cX^{\cf}_{\lsup LT}$
by one of its irreducible component
$\cY$, and we want to show 
$\cX_{\lsup LB}^{\cf}\underset{\cX_{\lsup LT}^{\cf}}{\times}\cY$
is irreducible.

Using the root space decomposition of $\Lie U$,
we can write it as a direct sum
$\Lie U = \oplus W_\alpha$,
where each $W_\alpha = \prod _{\gamma\in \Gal_{\bQ_p}}U_{\gamma\cdot \alpha}$ is the sum of the Galois orbit of a single root space
$U_\alpha$
such that each $H^0(\Gal_{\bQ_p}, W_\alpha)
\cong H^0(\Gal_L, U_\alpha)$
can only be either trivial or $1$-dimensional (by Shapiro's lemma).

Fix $\alpha\in \Delta(B, T)$.
There are three possibilities
\begin{enumerate}
\item $\dim H^0(\Gal_{\bQ_p}, W_\alpha)$ is always $0$ over $\cY$;
\item $\dim H^0(\Gal_{\bQ_p}, W_\alpha)$ is always $1$ over $\cY$, and
\item $\dim H^0(\Gal_{\bQ_p}, W_\alpha)$ vary over $\cY$.
\end{enumerate}
Denote by $\Delta_v$ the subset of simple roots
that lies in possibility (3).
We can thus divide $\cY$ into $2^{|\Delta_v|}$
locally closed substacks $\cY_Z, Z\subset \Delta_v$
where $\dim H^0(\Gal_{\bQ_p}, W_\alpha)=1$
if and only if $\alpha\in Z$.

The constancy of $\rk H^0(\Gal_{\bQ_p}, W_\alpha)$
implies
$\cX_Z:=\cX_{\lsup LB}^{\cf}\underset{\cX_{\lsup LT}^{\cf}}{\times}\cY_Z$
is Zariski-locally on $\cY_Z$ isomorphic to the product stack 
$\cY_Z \times [H^1(\Gal_{\bQ_p}, \Lie U)/
H^0(\Gal_{\bQ_p}, \Lie U)]$.
Indeed by induction on short exact sequences
$0\to \Lie U/W_\alpha\to\Lie U\to W_\alpha$,
it is harmless to replace $\Lie U$ by $W_\alpha$
and replace $\lsup LB$ by $W_\alpha\rtimes \lsup LT$.
Note that $H^i(\Gal_{\bQ_p}, W_\alpha)$
is a locally free sheaf over $\cY_Z$, $i=0,1,2$.
Then over a Zariski cover of $\cY_Z$
where $H^i(\Gal_{\bQ_p}, W_\alpha)$ is free,
$\cX_Z$ is a trivial $[H^1(\Gal_{\bQ_p}, W_\alpha)/
H^0(\Gal_{\bQ_p}, W_\alpha)]$-bundle.
For the same reason, each $\cX_Z$ is an irreducible stack.

By the Euler characteristic for Galois cohomology,
$$\dim [H^1(\Gal_{\bQ_p}, \Lie U)/H^0(\Gal_{\bQ_p}, \Lie U)]
= \dim H^2(\Gal_{\bQ_p}, \Lie U) + \dim U = \dim U$$
In particular, 
$\dim \cX_Z = \dim U + \dim \cY_Z
\begin{cases}
=\dim U & Z = \emptyset \\
<\dim U & Z \not= \emptyset
\end{cases}
$
by the irreducibility of $\cY$.
Since $\cX_{\lsup LB}^{\cf}$ is equidimensional,
its irreducibility follows from the
irreducibility of $\cX_\emptyset$.
\end{proof}

\begin{lem}
Each irreducible component of $\cX_{\lsup LB}^{\cf}$
has a dense substack which
embeds in $\cX^{\cf}$.
\label{lem:mns}
\end{lem}

\begin{proof}
The idea is to consider
$\cX_{\lsup LB}^{\cf}$
with $\cX_{\lsup LB \cap \lsup LB'}$ removed,
where $\lsup LB'\ne \lsup LB$ is another Borel subgroup.
By the Bruhat decomposition,
we only need to remove $\lsup LB'$
of the form $w\lsup LB'w^{-1}$
where $w$ is a Weyl group element.
Note that $\dim \cX_{\lsup LB \cap \lsup LB'} <\dim \cX_{\lsup LB}$ by Lemma \ref{lem:dim-UG}.
Let $\bar\rho\in \cX_{\lsup LB}^{\cf}-\cup_{B'}\cX_{\lsup LB \cap \lsup LB'}=:\mathcal{Z}$
and let $g\in \bG$ be an automorphism of $\bar\rho$.
We necessarily have $\bar\rho \in \cX_{\lsup LB \cap g\lsup LBg^{-1}}$,
which forces $g$ to be in the normalizer of $\brB$,
which is $\brB$ itself.
Therefore $\mathcal{Z}\to\cX_{\lsup LG}$
is a monomorphism.
\end{proof}

\begin{prop}
There is a natural bijection between
irreducible components of $\cX^{\cf}$
and $\cX_{\lsup LT}^{\cf}$.
\label{prop:LG-LT}
\end{prop}

\begin{proof}
Combine Lemma \ref{lem:mns}, Lemma \ref{lem:LB-LT} and
Lemma \ref{lem:dim-Xcf}.
\end{proof}

\begin{appendices}
\section{Algebraic characters and crystalline parameters}
\label{sec:A}

Let $T$ be a torus over $K$.
We also denote by $T$ its {\it standard}
integral model,
uniquely determined by $T(\cO_K)=T(K)^1$
(\cite[Fact B.4.3]{KP}).
The Iwahori model $T^\circ$ is the neutral component of $T$,
and we
denote by $\underline{T}$
the special fiber of the Iwahori
model $T^\circ$.
If $T$ is unramified, then $T=T^\circ$
and
$X^*(T)=X^*(\underline{T})\cong\bZ^{\rk T}$.

Fix an embedding $\bar K \to \bC$.
We say $T(K)\to \bC^{\times}$
is crystalline / semistable / de Rham
if the corresponding parameter
$\Gal_K\to\lsup LT(\bC)$
is crystalline / semistable / de Rham.

\begin{defn}
$f:T(K)\to \bC^\times$
is said to be an {\it algebraic character}
if there exists a cocharacter
$(\chi_\bullet)=(\chi_\sigma)_{\sigma: K \hookrightarrow \bC }
\in X^*(\Res_{K/\bQ_p}T)$
such that
$f(x) = \prod_{\sigma: K \hookrightarrow \bC}
\sigma(x^{\chi_\sigma})$
for all $x\in K$.

We will write $f=\ev(\chi_\bullet)$.
If $\chi\in X^*(T)$,
we also denote by $\chi$
the cocharacter $(\chi,0,0,\dots,0)\in X^*(\Res_{K/\bQ_pT})$
by abuse of notation.
\end{defn}

\begin{prop}
$f:T(K)\to \bC^\times$
is crystalline
if and only if $T(\cO_K) \to \bC^\times$
is an algebraic character in the sense
that there exists $\chi \in X^*(T)$
such that 
$f|_{T(\cO_K)}$ is the restriction
of an algebraic character
to $T(\cO_K)$.
\end{prop}

\begin{proof}
Since 
$T(\cO_K) \supset T(K)^0 = \underset{{E/K~\text{unramified}}}{\varinjlim}T(E)$,
crystallinity is completely determined
by the restriction of $f$ to $T(\cO_K)$;
algebraic characters $\ev(\chi_\bullet)$ are clearly
crystalline of Hodge type $\chi_\bullet$.
Conversely, the theory of weakly admissible
filtered $\varphi$-module together with
Dieudonne-Manin classification
imply
the Hodge type $\chi_\bullet$ completely determines
crystalline characters $T(\cO_K)\to \bC^\times$,
so $\ev(X^*(\Res_{K/\bQ_p}T))$ exhausts all
crystalline characters.
\end{proof}

\begin{prop}
Assume $\underline{T}$
splits over $\bF_{p^{nm}}$.
The following diagram is commutative:
\begin{equation*}
\begin{tikzcd}
0 \arrow[r] & X^*(\underline{T}) \arrow[r, "F-1"] 
\arrow[d, equal] & X^*(\underline{T}) 
\arrow[r] \arrow[d, "N", hook] & 
\operatorname{Hom}(T(\mathbb{F}_{p^n}), \bar{\mathbb{F}}^\times)
\arrow[d, "N", hook] 
\arrow[r] &0\\
0 \arrow[r] & X^*(\underline{T}) \arrow[r, "F^m-1"] & 
X^*(\underline{T}) 
\arrow[r, "\overline\ev"] & 
\operatorname{Hom}(T(\mathbb{F}_{p^{mn}}), \bar{\mathbb{F}}^\times)
\arrow[r] &0
\end{tikzcd}
\end{equation*}
Here $\overline\ev$ sends $\chi$ to $T(\bF_{p^{nm}})
\xrightarrow{x\mapsto x^\chi}
\Gm(\bF_{p^{nm}})\subset\bFp^{\times}$,
$N = 1+F+\dots+F^{m-1}$,
and $F=p\pi^{-1}$ is the absolute Frobenius structure
on $\underline{T}$.
All vertical arrows are injective.
\label{prop:dl}
\end{prop}

\begin{proof}
The short exactness of the bottom row is clear
because of the splitness of $\underline{T}$
over $\bF_{p^{nm}}$.
The rest of the proof follows from
\cite[Section 5.3]{DL}.
\end{proof}

\begin{lem}
$\Hom(T(K), \bar K^\times)\to \Hom(T(\cO_K), \bar K^\times)$
is surjective.
If $T$ is unramified, then 
$$\Hom(T(\cO_K), \bFp^\times)=
\Hom(\underline{T}(\bF_{p^n}), \bFp^\times)$$
where $\bF_{p^n}$ is the residue field of $\cO_K$.
\end{lem}
\begin{proof}
The first claim is because $\bar K^\times$
is a divisible group,
and $\Hom(-, \bar K^\times)$ is exact.
The second claim is clear.
\end{proof}

Let $E/K$ be a finite extension.
We say $T(K)\to \bC^\times$
{\it extends} $T(E)\to\bC^\times$
if, after passing to
    the Langlands dual side,
$\Gal_K\to \lsup LT(\bC)$
restricts to $\Gal_E\to\lsup LT(\bC)$.

\begin{lem}
Assume $T$ is unramified and splits over $K_m$.

(1)
Denote by $\Frob_K$ an arithmetic Frobenius element of $\Gal_K$.
Let $\rho:\Gal_{K_m}\to \bT(K_m)$
be a continuous homomorphism. 
Then $\rho(\Frob_K\circ \gamma \Frob_K^{-1})=\Frob_K(\rho(\gamma))$
for all $\gamma\in \Gal_{K_m}$.

(2)
Let
$\ev(\chi_\bullet):T(K_{m})\to K_m^{\times}$
be an algebraic character
defined by $\chi_\bullet\in X^*(\Res_{K/\bQ_p}T)$.
$\ev(\chi_\bullet)$ can be extended to
$T(K)\to K_m^\times$
if and only if $(\Frob_K \pi^{-1}-1)(\chi_\bullet)=0$.
\label{lem:alg-ext}

(3)
$(\Frob_K \pi^{-1} - 1)\ev(\chi_\bullet)\equiv
\ev((F-1)(\chi_\bullet))$
as a character in characteristic $p$.
\end{lem}

\begin{proof}
(1) It is \cite[Corollary 5.2.2]{L22}.

(2)
Pass to the dual side, and suppose we have a parameter
$\rho:\Gal_{K_m}\to \lsup LT(K_m)$.
$\rho$ can be extended to $\Gal_K$
if and only we can choose an element $X\in \bT(K_m)\rtimes \Frob_K$
such that after setting $\rho(\Frob_K)=X$,
we get $\rho(\Frob_K\circ \gamma \circ \Frob_K^{-1})=
X\rho(\gamma)X^{-1}$ ($*$)
for all $\gamma\in \Gal_{K_m}$.
Since $T$ splits over $K_m$, $\rho(\gamma)\in \bT(K_m)$
and therefore $X\ev(\chi_\bullet)X^{-1} = \ev(\pi(\chi_\bullet))$;
on the other hand,
$\ev(\chi_\bullet)\circ\Ad(\Frob_K)=
\Frob_K\ev(\chi_\bullet)
=\ev(\Frob_K(\chi_\bullet))$.
So ($*$) is equivalent to
$\Frob_K\pi^{-1}(\chi_\bullet)=(\chi_\bullet)$.

(3) Clear.
\end{proof}

\begin{thm}
\label{thm:crys-lift}
If $T$ is an unramified group,
then any character $\bar f:T(K)\to \bFp^\times$
admits a crystalline lift $T(K)\to \cO_\bC^\times$.
\end{thm}

\begin{proof}
Denote by $\bar\chi\in \Hom(T(\bF_{p^n}), \bFp^\times)$
the character that inflates to
of $\bar f|_{T(\cO_K)}$.
By Proposition \ref{prop:dl},
there exists a character $\chi\in X^*(\underline{T})$
such that $\bar\ev(N(\chi))=N(\bar f)$.
Set
$$
N(\ev(\chi))
:=(1+\Frob_K\pi^{-1} +\dots+\Frob_K^m\pi^{-m})\ev(\chi),
$$
which defines an algebraic character
$T(K_m) \to K_m^\times$
lifting $\bar\ev(N(\chi))$.
We have 
$$
(\Frob_K\pi^{-1}-1)N(\ev(\chi))
=(\Frob_K^m\pi^{-m}-1)\ev(\chi)=0\ev(\chi)=0.
$$
So by Lemma \ref{lem:alg-ext},
it extends to $g:T(K)\to K_m^\times$.
By the injectivity part of Proposition \ref{prop:dl},
$\bar g:T(K)\to \cO_{K_m}^\times\to \bF_{p^n}^\times$
and $\bar f$ are equal over $T(\cO_K)$.
Finally, denote by $[\bar f \bar g^{-1}]:T(K)\to K_m^\times$
the Teichm\"uller lift of $(\bar f \bar g^{-1})$,
the product $[\bar f \bar g^{-1}] g$
is a crystalline lift of $\bar f$.
\end{proof}

\begin{lem}
Set $K=\bQ_p$.
Let $\lambda\in X^*(T)$
be an arbitrary character.
There exists a potentially crystalline
character $T(\bQ_p)\to \cO_{\bC}^\times$
of Hodge type $\lambda$.
\label{lem:exist-de-rham}
\end{lem}

\begin{proof}
Use the embedding $T\hookrightarrow\Res_{L/\bQ_p}T$,
and the character
$\ev(\lambda):\Res_{L/\bQ_p}T(\bQ_p)=T(L)\to \cO_{\bC}^\times$.
\end{proof}

\section{Moduli of de Rham parameters: the reductive case}

\label{sec:moduli-ss}
Let $G$ be a quasi-split reductive group over $K$
that splits over $E$,
where $E$ is a tamely ramified extension of $K$.
Write $\lsup LG=\wh G\rtimes \Gal(E/K)$
for the Langlands dual group of $G$,
where $\wh G$ is the pinned dual group of $G$ defined over $\spec \bZ$.

We freely use the moduli stacks constructed in \cite{L23B}, see \cite[Table 1.7.1]{L23B}.

\subsection{Breuil-Kisin modules with $\wh G$-structure}~

Let $\cO\supset \cO_K$ be a DVR over $\bZ_p$.

Denote by $k$ the residue field of $K$.
Let $A$ be a $\bZ_p$-algebra.
For each choice of a compatible family
$\pi^{1/p^\infty} = (\pi^{1/p^n})_{n\in \bZ_+}$
of $p$-power roots of a uniformizer of $K$
in $\bQp$,
we define an embedding
\begin{align*}
(W(k)\otimes_{\bZ_p} A)[[u]] & \to \bA_{\Inf, A} \\
u&\mapsto [\pi^\flat]
\end{align*}
where $\pi^\flat = \invlim{n}\pi^{1/p^n} \in \cO_{\bC}^\flat$.
Denote by $\fS_{\pi^\flat, A}$
the image of the embedding above.

\subsubsection{Definition}
A projective {\it Breuil-Kisin module} with $A$-coefficients
is a finitely generated projective 
$\fS_{\pi^\flat, A}$-module $\fM$,
equipped with
a $\varphi$-semi-linear morphism
$\phi_{\fM}:\fM\to \fM$
such that $1\otimes \phi_\fM:\varphi^*\fM[1/E_{\pi^\flat}]\to \fM[1/E_{\pi^\flat}]$
is a bijection.
Here $E_{\pi^\flat}$ is the Eisenstein polynomial corresponding
to $\pi$.

We say $(\fM, \phi_\fM)$ is an {\it effective} Breuil-Kisin module
if $\phi_\fM(\fM)\subset \fM$.
We say a Breuil-Kisin module $(\fM, \phi_\fM)$
has {\it height $h$}
if $${E_{\pi^\flat}^h}\fM\subset\Img(1\otimes \fM)\subset \frac{1}{E_{\pi^\flat}^h}\fM.$$

Following the notation of \cite{EG23},
we denote by $\cC_{\pi^\flat, d, h}$
the (limit-preserving) moduli stack of effective Breuil-Kisin modules
of height at most $h$
for the uniformizer $\pi$
(see \cite[4.5.7]{EG23} for the definition).
Also denote by $\cC_{\pi^\flat, d, h}^a$
the base-changed stack
$\cC_{\pi^\flat, d, h} \otimes_{\bZ_p}\cO/\varpi^a$.

Denote by $\clR_{\pi^\flat, d}$ the moduli of rank-$d$ \'etale $\varphi$-modules
for the $\varphi$-ring $\fS_{\pi^\flat}[1/u]$,
and denote by $\clR_{\pi^\flat, d}^a$ the base-changed version of $\clR_{\pi^\flat, d}$.

In \cite[Section 3]{L23B}, we constructed
the moduli stack
$\cC_{\pi^\flat, \wh G, h}$
of Breuil-Kisin modules with $\wh G$-structure,
and the moduli stack
$\clR_{\pi^\flat, \wh G}$
of rank-$d$ \'etale $\varphi$-modules with $\wh G$-structure.
Note that when $\wh G=\GL_d$,
there is a canonical monomorphism
$$
\cC_{\pi^\flat, d, h}
\hookrightarrow
\cC_{\pi^\flat, \GL_d, h}
$$
as $\cC_{\pi^\flat, \GL_d, h}$
classifies Breuil-Kisin modules
of height $h$ that are not necessarily effective.
Moreover, the $h$-th Tate twist morphism
$$
\cC_{\pi^\flat, \GL_d, h} \xrightarrow{\cong} \cC_{\pi^\flat, d, 2h}
$$
sending $F$ to $F\times^{\GL_d}\Ga^{\oplus d}(h)$
is an equivalence of $2$-categories.

\subsubsection{Canonical extensions of $\Gal_{K_\infty}$-actions}
\label{par:can}
Write $K_s$ for $K(\pi^{1/p^s})$,
$0\le s\le \infty$.
We also write $K_{\pi^\flat, s}=K_s$
to emphasize the choice of $\pi^\flat$.

By Fontaine's theory (\cite[Proposition 2.7.8]{EG23}), we can attach to 
a Breuil-Kisin module a $(\varphi, \Gal_{K_\infty})$-module.
By \cite[Proposition 4.5.8]{EG23},
the $\Gal_{K_\infty}$-action
on a $(\varphi, \Gal_{K_\infty})$-module admits
a canonical extension to $\Gal_{K_s}$,
where $s$ is any integer greater
than the constant $s(a, h, N)$ defined
in \cite[Lemma 4.3.3]{EG23}
(where $N>\frac{e(a+h)}{p-a}$ is a fixed number).
More precisely, there exists a commutative diagram
\begin{equation}
\label{eqn:can-d}
\xymatrix{
& \cC^a_{\pi^\flat, d, h} \ar[dl]\ar[d]\\
\cX_{K_s, d}^a \ar[r] & \clR_{\pi^\flat, d}^a
}
\end{equation}

where $\cX^a_{K_s,d}=\cX_{K_s, d}\otimes_{\bZ_p}\cO/\varpi^a$ is the base-changed Emerton-Gee stack.

\subsubsection{Proposition}
For each $s>s(a, 2h, N)$,
there exists a commutative diagram
$$
\xymatrix{
& \cC^a_{\pi^\flat, d, h} \ar@{^{(}->}[d]\\
& \cC^a_{\pi^\flat, \GL_d, h} \ar[dl]\ar[d]\\
\cX_{K_s, d}^a \ar[r] & \clR_{\pi^\flat, d}^a
}
$$
extending Diagram (\ref{eqn:can-d}).

\begin{proof}
It follows immediately from the discussion before
Paragraph \ref{par:can}.
\end{proof}

\subsubsection{Definition}
Fix once for all an embedding
$\wh G\hookrightarrow \GL_d$.
We say a Breuil-Kisin module $(F, \phi_F)$
with $\wh G$-structure
has {\it height $h$}
if $(F, \phi_F)\times^{\wh G}\GL_d$
has height $h$.

\subsubsection{Lemma}
For each $s\ge 0$,
the morphism
$$
\text{($\dagger$)}\hspace{10mm}\cX_{K_s, \wh G}^a
\hookrightarrow
\cX_{K_s, \GL_d}^a
\underset{\clR_{\pi^\flat, \GL_d}^a}{\times}\clR_{\pi^\flat, \wh G}^a
$$
is a closed immersion.

\begin{proof}
Write $V=\Ga^{\oplus d}$ for the trivial vector space
of rank $d$ so that $\GL_d=\GL(V)$.
For ease of notation, set $s=0$.
Fix a finite type $\cO/\varpi^a$-algebra
$A$.
$\cX_{K, \GL_d}^a(A)$
is equivalent to the groupoid
of projective \'etale
$(\varphi, \Gal_{K})$-modules
of rank $d$ with $A$-coefficients,
and $\clR_{\pi^\flat, \GL_d}^a(A)$
is equivalent to the groupoid
of projective \'etale
$(\varphi, \Gal_{K_\infty})$-modules
of rank $d$ with $A$-coefficients.
An object of 
$\cX_{K, \GL_d}^a
\underset{\clR_{\pi^\flat, \GL_d}^a}{\times}\clR_{\pi^\flat, \wh G}^a(A)$
is a tuple
$$
(M, F, \iota):=((M, \phi_M, \rho_M), (F, \phi_F, \rho_{F,\infty}),
\iota: (M, \phi_M, \rho_{M}|_{\Gal_{K_\infty}})\cong
(F, \phi_F, \rho_{F, \infty})\times^{\wh G}V),
$$
where $(M, \phi_M, \rho_M)\in \cX_{K, \GL_d}^a(A)$
and $(F, \phi_F, \rho_{F,\infty})\in \clR_{\pi^\flat, \wh G}^a(A)$.
A morphism
$$
(M_1, F_1, \iota_1)\to (M_2, F_2, \iota_2)
$$
is a pair $(g: M_1\to M_2, f:F_1\to F_2)$ such that
$$
g = \iota_2^{-1}\circ (f\times^{\wh G}V)\circ \iota_1.
$$
The morphism  ($\dagger$) can be explicitly written as
$$
(F, \phi_F, \rho_F)
\mapsto (F\times^{\wh G}V, F, \id).
$$
From the description above, the morphism ($\dagger$)
is clearly faithful;
indeed ($\dagger$) is fully faithful:
a $\wh G$-torsor morphism
$F_1\to F_2$ respects the $\Gal_K$-action
if and only if $F_1\times^{\wh G}\GL_d
\to F_1\times^{\wh G}\GL_d$
respects the $\Gal_K$-action
since $F_i\hookrightarrow F_i\times^{\wh G}\GL_d$
is a closed subscheme ($i=1, 2$).

By \cite[Lemma 10.3.2]{L23B}, $(\dagger)$
is of strong Ind-finite type
in the sense of \cite[Definition 10.3.1]{L23B}.
To show $(\dagger)$ is a closed immersion,
it suffices to show it is proper using the (Noetherian) valuative criterion.
Let $\Lambda$ be a discrete valuation ring
over $\bFp$ with fraction field $\Omega$.
By \cite[Tag 0ARL]{Stacks},
it is harmless to assume $\Lambda$ is a complete discrete valuation ring;
thus $\fS_{\pi^\flat,\Omega}$
is a disjoint union of the spectrum of Noetherian complete regular local rings over $\bFp$.
By the Grothendieck-Serre conjecture (see the main theorem of \cite{FP15}),
all $\wh G$-torsors over $\fS_{\pi^\flat,\Lambda}$
are trivial $\wh G$-torsors.
The valuative criterion can be checked by noticing
the simple fact that
$\wh G(\fS_{\pi^\flat,\Omega})\cap \GL(\fS_{\pi^\flat,\Lambda}) = \wh G(\fS_{\pi^\flat,\Lambda})$.
\end{proof}

\subsubsection{Lemma}
The diagram
$$
\xymatrix{
\cX_{K_s, \wh G}^a
\underset{\cX_{K_s, \GL_d}^a}{\times}\cC_{\pi^\flat, \GL_d, h}^a\ar@{^{(}->}[r]\ar[d]
 &
\clR_{\pi^\flat, \wh G}^a
\underset{\clR_{\pi^\flat, \GL_d}^a}{\times}\cC_{\pi^\flat, \GL_d, h}^a \ar[d]\\
\cX_{K_s, \wh G}^a\ar@{^{(}->}[r]
 &
\clR_{\pi^\flat, \wh G}^a
\underset{\clR_{\pi^\flat, \GL_d}^a}{\times}
\cX_{K_s, \GL_d}^a
}
$$
is Cartesian.

\begin{proof}
Note that
$
A\underset{B\underset{C}{\otimes}D}{\otimes}(B\underset{C}{\otimes}E)
=A\underset{B\underset{C}{\otimes}D}{\otimes}(B\underset{C}{\otimes}D)\otimes_DE
=A\otimes_DE.
$
\end{proof}

\subsubsection{Canonical extensions of $\Gal_{K_\infty}$-actions, the $\wh G$-version}~

By \cite[Corollary 3.7.2]{L23B},
the morphism
$$
\cC_{\pi^\flat, \wh G, h}
\hookrightarrow
\cC_{\pi^\flat, \GL_d, h}
\underset{\clR_{\pi^\flat, \GL_d}}{\times}\clR_{\pi^\flat, \wh G}
$$
is a closed immersion.
Define the stack
$\cC_{\pi^\flat, K_s, \wh G, h}^a$
so that the following diagram
is Cartesian
\begin{equation}
\label{eq:def-CpiKGh}
\xymatrix{
\cC_{\pi^\flat, K_s, \wh G, h}^a\ar@{^{(}->}[rr]
\ar@{^{(}->}[d]
&&  \cC_{\pi^\flat, \wh G, h}^a
\ar@{^{(}->}[d]
\\
\cX_{K_s, \wh G}^a
\underset{\cX_{K_s, \GL_d}^a}{\times}\cC_{\pi^\flat, \GL_d, h}^a
\ar@{^{(}->}[rr]
&&
\clR_{\pi^\flat, \wh G}^a
\underset{\clR_{\pi^\flat, \GL_d}^a}{\times}\cC_{\pi^\flat, \GL_d, h}^a
}
\end{equation}
Note that all arrows in the diagram above are closed immersions.
The diagram above can be interpreted as follows:
for all Breuil-Kisin module with $\wh G$-structure,
the $\Gal_{K_{\infty}}$-action
can be extended to a $\Gal_{K_s}$-action canonically
in a way not necessarily compatible with the
$\wh G$-structure;
and the condition that the canonical extension
is compatible with the $\wh G$-structure is a closed condition.

Define the stack
$\cC_{\pi^\flat, s, \wh G, h}$
so that the following diagram
is Cartesian
$$
\xymatrix{
\cC_{\pi^\flat, s, \wh G, h}^a\ar[rr]
\ar[d]
&&
\cC_{\pi^\flat, K_s, \wh G, h}^a
\ar[d]
\\
\cX_{K, \wh G}^a
\ar[rr]
&&
\cX_{K_s, \wh G}^a
}
$$
is Cartesian.
The stack
$\cC_{\pi^\flat, s, \wh G, h}^a$
can be interpreted as the moduli stack
of Breuil-Kisin modules $(\fM, \phi_{\fM})$
together with an enhancement: a $(\varphi, \Gal_{K_s})$-module
$(\fM\underset{\fS_{\pi^\flat, A}[1/u]}{\otimes}W(\bF^\flat)_A, \phi_{\fM}\otimes 1, \rho)$
such that $\rho\times^{\wh G}{\GL_d}|_{\Gal_{K_s}}$
is the canonical $\Gal_{K_s}$-action.

\subsubsection{Lemma}
\label{lem:cCsGh-diag}
(1)
The morphism
$\cC^a_{\pi^\flat, K_s, \wh G, h}
\to \cX_{K_s, \wh G}^a$
is representable by algebraic spaces,
proper, and of finite presentation.

(2)
The diagonal of $\cC_{\pi^\flat, s, \wh G, h}^a$
is affine and of finite presentation.

\begin{proof}
(1) It follows from the diagram (\ref{eq:def-CpiKGh})
and \cite[Lemma 4.5.9]{EG23}.

(2) It follows from part (1),
\cite[Theorem 7.1.2]{L23B},
and \cite[Lemma 4.5.14]{EG23}.
\end{proof}

\subsection{Breuil-Kisin-Fargues modules with $\wh G$-structure}~

Let $L/\Qp$ be a finite Galois extension,
Let $\cO\supset \cO_L$ be a DVR over $\bZ_p$.

Denote by $l$ the residue field of $L$.
Let $A$ be a $\bZ_p$-algebra.
For each choice of a compatible family
$\pi^{1/p^\infty} = (\pi^{1/p^n})_{n\in \bZ_+}$
of $p$-power roots of a uniformizer of $L$
in $\bQp$,
we define an embedding
\begin{align*}
(W(l)\otimes A)[[u]] & \to \bA_{\Inf, A} \\
u&\mapsto [\pi^\flat]
\end{align*}
where $\pi^\flat = \invlim{n}\pi^{1/p^n} \in \cO_{\bC}^\flat$.
Denote by $\fS_{\pi^\flat, A}$
the image of the embedding above.

\subsubsection{Definition}
A  projective {\it Breuil-Kisin-Fargues} module
with $A$-coefficients
is a finitely generated projective
$\bA_{\Inf, A}$-module
$\fM^\Inf$, equipped with
a $\varphi$-semi-linear endomorphism
$\phi_{\fM^\Inf}:\fM^\Inf\to \fM^\Inf$ such that
$1\otimes \phi_{\fM^\Inf}:\varphi^*\fM^\Inf[1/\xi]\to \fM^\Inf[1/\xi]$
is a bijection.
Here $\xi$ is a generator of $\ker(\theta:\bA_\Inf \to \cO_{\bC})$
(\cite[Definition 4.22]{BMS18}).

We say a Breuil-Kisin-Fargues module $\fM^\Inf$
{\it descends to $\fS_{\pi^\flat, A}$}
if there is a Breuil-Kisin module
$\fM_{\pi^\flat}\subset(\fM^\Inf)^{\Gal_{K_{\pi^\flat,\infty}}}$
such that
$\bA_{\Inf, A}\underset{\fS_{\pi^\flat, A}}{\otimes}\fM_{\pi^\flat}
=\fM^\Inf$.
We say $\fM^\Inf$ {\it admits all descents over $L$}
if it descends to $\fS_{\pi^\flat, A}$
for every choice of $\pi^\flat$ (for every choice of $\pi$)
and if furthermore
\begin{itemize}
\item 
[(1)]
the $W(l)\otimes_{\bZ_p}A$-submodule
$\fM_{\pi^\flat}/[\pi^\flat]\fM_{\pi^\flat}$
of $(W(\bar l)\otimes A)\otimes_{\bA_{\Inf, A}}\fM^\Inf$
is independent of the choice of $\pi$
and $\pi^\flat$;
\item[(2)]
the $\cO_L\otimes_{\bZ_p}A$-submodule
$\varphi^*\fM_{\pi^\flat}/E_{\pi^\flat}\varphi^*\fM_{\pi^\flat}$
of $\cO_{\bC, A}\otimes_{\bA_{\Inf, A}}\varphi^*\fM^\Inf$
is independent of the choice of $\pi$ and $\pi^\flat$.
\end{itemize}

\subsubsection{Lemma}
\label{lem:BKF-des-tan}
The category of Breuil-Kisin-Fargues modules
with $A$-coefficients that
admits all descents over $L$
is an exact, rigid, symmetric monoidal category.

\begin{proof}
Note that the property of
admitting all descents over $L$ is preserved
under tensor products and duals.
\end{proof}

\subsubsection{Definition}
Let $A$ be a $p$-adically complete $\cO$-algebra
which is topologically of finite type.
A {\it Breuil-Kisin-Fargues $\Gal_L$-module
with $A$-coefficients}
is a Breuil-Kisin-Fargues module $(\fM^\Inf, \phi_{\fM^\Inf})$
with $A$-coefficients equipped with
a continuous semilinear $\Gal_L$-action
that commutes with $\phi_{\fM^\Inf}$.

\subsubsection{Connection to semistable Galois representations}~

Let $M$ be an \'etale $(\varphi, \Gal_L)$-module.
By \cite[Section 2.7]{EG23},
we can attach a $\Gal_L$-representation $V(M)$ to $M$.
By \cite[Theorem F.11]{EG23},
if $V(M)$ is semistable with Hodge-Tate weights
    in $[0,h]$,
    then there exists a unique Breuil-Kisin-Fargues $\Gal_L$-module $\fM^\Inf$
    of height at most $h$
    that admits all descents over $L$
    such that
    $\fM^\Inf\otimes_{\bA_\Inf}W(\bC^\flat)=M$.
Moreover, by \cite[Proposition 4.4.1]{EG23},
there exists a constant $s'(L, a, h, N)$
(where $N>\frac{e (a+h)}{p-a}$ is a fixed constant),
such that for any choice $\pi^\flat$
with corresponding descent $\fM_{\pi^\flat}$
and for any $s> s'(L, a, h, N)$,
the restriction to $\Gal_{L_{\pi^\flat, s}}$
of the action of $\Gal_{L}$ on $\fM^\Inf \otimes_{\cO}\cO/\varpi^a$ agrees with the canonical action
considered in Paragraph \ref{par:can}.

\subsubsection{Definition}
A {\it Breuil-Kisin-Fargues $\Gal_L$-module
with $A$-coefficients with $\wh G$-structure}
is a faithful, exact, symmetric monoidal functor
from $\frep_{\wh G}$ to the category of
Breuil-Kisin-Fargues $\Gal_L$-module
with $A$-coefficients.

\subsubsection{Proposition}
\label{prop:ss-BKF}
Let $F$ be an \'etale $(\varphi, \Gal_L)$-module
with $\wh G$-structure.
Then $V(F)$ is a semistable $\wh G$-valued Galois representation
of $\Gal_L$ such that $V(F)\times^{\wh G}\GL_d$
has Hodge-Tate weights in $[-h, h]$
if and only if there exists a (necessarily) unique
Breuil-Kisin-Fargues $\Gal_L$-module $F^\Inf$
with $\wh G$-structure
which admits all descents over $L$ and which satisfies
$F=F^\Inf \otimes_{\bA_{\Inf}}W(\bC^\flat)$,
and such that
$F^\Inf\times^{\wh G}\GL_d$ is of height at most $h$.

Moreover, 
for any choice $\pi^\flat$
with corresponding descent $F_{\pi^\flat}$
and for any $s> s'(L, a, 2h, N)$,
the restriction to $\Gal_{L_{\pi^\flat, s}}$
of the action of $\Gal_{L}$ on $(F^\Inf\times^{\wh G}\GL_d) \otimes_{\cO}\cO/\varpi^a$ agrees with the canonical action
considered in Paragraph \ref{par:can}.

\begin{proof}
Let $x\in \frep_{\wh G}$.
By \cite[Proposition 5.3.2, Definition 5.3.1]{Lev13},
$V(F)\times^{\wh G}x=V(F\times^{\wh G}x)$
is semistable.
Fix a choice of $\pi^\flat$.
By \cite[Theorem F.11]{EG23},
$F\times^{\wh G}x$
corresponds to a unique Breuil-Kisin-Fargues $\Gal_L$-module
$F^\Inf_x$.
The uniqueness implies the association
$F^\Inf:x\mapsto F^\Inf_x$ is functorial and monoidal.
Since $F^\Inf \otimes_{A_\Inf}W(\bC)=F$
is faithful and exact
and that $A_\Inf$ is a subring of $W(\bC)$,
$F^\Inf$ is automatically faithful and lef-exact.
It remains to show $F^\Inf$ is right-exact,
and we do it by the bundle extension technique.

By \cite[Lemma 4.2.8]{EG23},
$F^\Inf_x$ descends for each $\pi^\flat$ uniquely 
to a Breuil-Kisin module $F^{\pi^\flat}_x$.
Since $\fS_{\pi^\flat}\to A_\Inf$
is faithfully flat,
it suffices to show
the functor $F^{\pi^\flat}: x\mapsto F^{\pi^\flat}_x$
is right-exact.
By \cite[Lemma 4.2.22]{Lev13},
the functor $F^{\pi^\flat}[1/p]$
is exact, and thus
defines a $\varphi$-module with $\wh G$-structure $F_1$
over $\spec \fS_{\pi^\flat}[1/p]$.
The functor $F^{\pi^\flat}[1/u]$
is also exact
by almost \'etale descent
($\fS[1/u]\to W(\bC)$ is faithfully flat, see \cite[Proposition 2.2.14]{EG23}),
and thus defines a $\varphi$-module with $\wh G$-structure $F_2$
over $\spec \fS_{\pi^\flat}[1/u]$.
The two $\varphi$-modules $F_1$ and $F_2$
can be glued along the intersection of
$\spec \fS_{\pi^\flat}[1/p]$
and
$\spec \fS_{\pi^\flat}[1/u]$.
By \cite[Lemma 5.1.1]{Lev13} (see also the main result of \cite{An22}),
there exists unique extension $F_3$ of $F_1$ and $F_2$
to $\spec \fS_{\pi^\flat}$.
Since $F_3$, when regarded as a monoidal functor, is forced to be $F^{\pi^\flat}$,
we have finished showing $F^\Inf$
is a faithful, exact, symmetric monoidal functor.

The ``moreover'' part is \cite[Proposition 4.4.1]{EG23}.
\end{proof}

The proposition above motivates the following definition.

\subsubsection{Definition}
For each $h\ge 0$ and each finite Galois extension $L$ of $E$,
define $\cC_{L, \wh G, \sS, h}^a$
to be the limit-preserving stack over $\cO/\varpi^a$
such that for each finite type $\cO/\varpi^a$-algebra,
$\cC_{L, \wh G, \sS, h}^a(A)$
classifies Breuil-Kisin-Fargues $\Gal_L$-modules $F^\Inf$
with $A$-coefficients and $\wh G$-structure
which admits all descents over $L$
such that $F^\Inf \times^{\wh G}\GL_d$
is of height $h$ and is equipped with canonical $\Gal_L$-action in the sense of Proposition \ref{prop:ss-BKF}.

\subsubsection{Lemma}
\label{lem:CLGssh}
Let $\pi_L^\flat$ be a compatible
system of $p$-power roots of a uniformizer
$\pi_L$ of $L$.
The morphism
\begin{equation}
\label{eq:CLGssh}
\cC_{L, \wh G, \sS, h}^a
\hookrightarrow
\cC_{L, \GL_d, \sS, h}^a
\underset{\cC_{\pi^\flat_L, s, \GL_d, h}^a}{\times}
\cC_{\pi^\flat, s, \wh G, h}^a
\end{equation}
is a closed immersion.

\begin{proof}
The morphism $\cC_{L, \wh G, \sS, h}^a\to\cC_{\pi^\flat, s, \wh G, h}^a$ is monomorphism by the definition
of $\cC_{\pi^\flat, s, \wh G, h}^a$.
By \cite[Proposition 4.5.17]{EG23},
(\ref{eq:CLGssh}) is a monomorphism.
Let $A$ be a finite type $\cO/\varpi^a$-algebra.
By Lemma \ref{lem:BKF-des-tan},
An object $F\in \cC_{\pi^\flat, s, \wh G, h}^a(A)$
lies in the essential image of $\cC_{L, \wh G, \sS, h}^a(A)$ if and only if for any $x\in \frep_{\wh G}$,
$F\times^{\wh G}x$ admits all descents over $L$ when regarded as a Breuil-Kisin-Fargues $\Gal_L$-module with $A$-coefficients.
By the proof of \cite[Proposition 4.5.17]{EG23},
the admitting all descents over $L$ is a closed condition.
\end{proof}

\subsubsection{Lemma}
\label{lem:BKF-1}
$\cC_{L, \wh G, \sS, h}^a$
is an algebraic stack of finite presentation over $\cO/\varpi^a$, and have affine diagonals.

\begin{proof}
It follows from Lemma \ref{lem:cCsGh-diag}
and Lemma \ref{lem:CLGssh}.
\end{proof}

\subsubsection{Lemma}
\label{lem:BKF-2}
$\cC_{L, \wh G, \sS, h} := \dirlim{a}~\cC_{L, \wh G, \sS, h}^a$
is a $p$-adic formal algebraic stack of finite presentation,
and have affine diagonals.
The morphism $\cC_{L, \wh G, \sS, h}\to \cX_{L, \wh G}$
is representable by algebraic spaces, proper
and of finite presentation.

\begin{proof}
It is a generalization of \cite[Theorem 4.5.20]{EG23},
and follows from Lemma \ref{lem:BKF-1}
and \cite[Proposition 4.5.17]{EG23}.
\end{proof}

\subsection{Breuil-Kisin-Fargues lattices with $\lsup LG$-structure}~

Now we generalize the previous discussions to
the case of non-split groups.
Fix a uniformizer $\pi=\pi_K$ of $K$.

\subsubsection{Lemma}
If $E$ is a tamely ramified Galois extension of $K$
of ramification index $e$,
then $\pi_E:=\pi_K^{1/e}$ is a uniformizer of $E$.

\begin{proof}
We first show that $\pi_E\in E'$ for some unramified extension $E'$ of $E$.
Let $\Pi\in \cO_E$ be an arbitrary uniformizer.
We have $\Pi^e= c \pi$ for some $c\in \cO_E^\times$.
Write $\bar c\in \kappa_E$ for the image
of $c$ in the residue field $\kappa_E$
and write $[\bar c]$ for the Teichm\"uller lift of $\bar c$.
Since $[\bar c]$ admits an $e$-th root $\mu$ in an unramified extension
of $E$, by replacing $\Pi$ by $\Pi \mu$, we can assume $c = 1+\pi x$ for some $x\in \cO_E$.
Thus $\frac{\Pi^e}{\pi_E^e}-1\in \Pi \cO_E$.
After multiplying $\Pi$ by an $e$-th root of unity,
we have $\frac{\Pi}{\pi_E}-1\in \Pi \cO_E$,
and thus $|\Pi - \pi_E| < |\Pi|$.
By Krasner's lemma (see, for example, \cite[Lemma 4.1.1]{L23B}),
we have $\pi_E\in E'[\Pi]=E'$.

Since $x^e=\pi$ is an Eisenstein polynomial, $K[\pi_E]$
is totally ramified over $K$. Therefore $\pi_E\in E$.
\end{proof}

Fix a tame Galois extension $E$ of $K$ of ramification index $e$.
By the lemma above, $\pi_E:=\pi^{1/e}$ is a uniformizer of $E$.
Let $\pi_E^\flat$ be a compatible system of $p$-power roots of $\pi_E$.
Set $\pi^\flat:=(\pi_E^\flat)^e$,
which is a compatible system of $p$-power roots of $\pi$.
Set $u:=[\pi^\flat]$ and $u_E:=[\pi_E^\flat]$.
We have
\begin{eqnarray*}
\fS_K:=&\fS_{\pi^\flat} &= \kappa_K[[u]] \\
\fS_E:=&\fS_{\pi_E^\flat} &= \kappa_E[[u_E]].
\end{eqnarray*}
It is clear that $\spec \fS_E[1/u_E]\to \spec \fS_K[1/u]$
is a Galois cover with Galois group
canonically identified with $\Gal(E/K)$.

\subsubsection{Semistable Galois representations and Breuil-Kisin-Fargues lattices}~

There exists an $\bA_\Inf$-linear isomorphism
\begin{eqnarray*}
A_\Inf \underset{\fS_K}{\otimes}\fS_E
&\to \prod_{\sigma\in \Gal(E/K)} A_\Inf\\
a\otimes b &\mapsto (a\sigma(b))_\sigma
\end{eqnarray*}
The diagonal map
$\Delta:A_\Inf \underset{\fS_K}{\otimes}\fS_E
\to A_\Inf$, $a\otimes b\mapsto ab$
is the unique $\Gal_E$-equivariant $A_\Inf$-linear
homomorphism such that
the composition
$$
A_\Inf \xrightarrow{a\mapsto a\otimes 1}
A_\Inf \underset{\fS_K}{\otimes}\fS_E
\xrightarrow{\Delta}
A_\Inf
$$
is the identity map.
Recall that $A_\Inf = W(\cO_{\bC^\flat})$;
we also need to consider the base-changed version
$$
W(\bC^\flat) \xrightarrow{a\mapsto a\otimes 1}
W(\bC^\flat) \underset{\fS_K}{\otimes}\fS_E
=
W(\bC^\flat)\underset{\bA_K}{\otimes}\bA_E
\xrightarrow{\Delta}
W(\bC^\flat).
$$

Let $(F, \phi_F, \rho_F)$
be an \'etale $(\varphi, \Gal_K)$-module
with $\lsup LG$-structure over $W(\bC^\flat)$
that corresponds to an $L$-parameter.
By \cite[Subsection 1.2]{L23B},
we have a $\Gal(E/K)$-equivariant isomorphism
\begin{equation}
\label{eq:wh-G-lvl}
\bar F:=F \times^{\lsup LG} \Gal(E/K)
\cong \spec W(\bC^\flat)\otimes_{\bA_K}\bA_E
\end{equation}
which we fix once for all.
Moreover $F\to \bar F$ defines a $\wh G$-torsor over $\bar F$.
Consider the pullback diagram
$$
\xymatrix{
F_\Delta \ar[d]\ar@{^{(}->}[r] & 
F\ar[d]
\\
\spec W(\bC^\flat) \ar@{^{(}->}[r]^{\hspace{6mm}\Delta} &
\bar F
}
$$
Since $\Delta$ is an $\Gal_E$-equivariant embedding,
$F_\Delta\hookrightarrow F$
is both $\varphi$-equivariant and $\Gal_E$-equivariant;
therefore $F_\Delta$ inherits
a structure of \'etale $(\varphi, \Gal_E)$-module
with $\wh G$-structure,
which we write
as $(F_\Delta, \phi_\Delta, \rho_\Delta)$
for simplicity.

Now assume 
the \'etale $(\varphi, \Gal_K)$-module
$(F, \phi_F, \rho_F)$
with $\lsup LG$-structure 
is potentially semistable;
suppose it becomes semistable
after restricting to $\Gal_L$
for some finite Galois extension $L$ of $E$.
By Proposition \ref{prop:ss-BKF},
the semistable \'etale $(\varphi, \Gal_L)$-module
$(F_\Delta, \phi_\Delta, \rho_\Delta|_{\Gal_L})$
admits a unique Breuil-Kisin-Fargues $\Gal_L$-lattice
$F_\Delta^\Inf$ with $\wh G$-structure.
Note that
$F_\Delta^\Inf$
is $\Gal_E$-invariant
(c.f. \cite[Corollary F.23]{EG23}),
and $F_\Delta^\Inf$ is indeed a 
Breuil-Kisin-Fargues $\Gal_E$-module with $\wh G$-structure.
Consider the following pushout diagram
$$
\xymatrix{
F_\Delta \ar@{^{(}->}[r]\ar[d] & 
F \ar@{=}[r] & 
\underset{\sigma\in \Gal(E/K)}{\coprod}\sigma(F_\Delta)
\ar[d]
\\
F_\Delta^\Inf \ar@{^{(}->}[rr]
&&
\underset{\sigma\in \Gal(E/K)}{\coprod}\sigma(F_\Delta^\Inf)
}
$$
Set $F^\Inf:=\underset{\sigma\in \Gal(E/K)}{\coprod}\sigma(F_\Delta^\Inf)$.
Concretely,
$\cO_{F^\Inf_\Delta}\subset \cO_{F_\Delta}$
(the structure sheaf)
is an $\bA_\Inf$-submodule,
and $\cO_{F^\Inf}$
is by definition the (direct) sum of the $\Gal(E/K)$-translations
of $\cO_{F_\Delta^\Inf}$
in $\cO_F$
(note that $F_\Delta$ is a connected component of $F$
and $\cO_{F_\Delta}$ is a direct summand of $\cO_F$).

\subsubsection{Proposition}
\label{prop:ss-BKF}
Let $F$ be an \'etale $(\varphi, \Gal_K)$-module
with $\lsup LG$-structure.
Assume $V(F)$ is a potentially semistable $\lsup LG$-valued $L$-parameter.
Then there exists a unique
Breuil-Kisin-Fargues $\Gal_K$-module $F^\Inf$
with $\lsup LG$-structure
which admits all descents over $L$ and which satisfies
$F=F^\Inf \otimes_{\bA_{\Inf}}W(\bC^\flat)$.

\begin{proof}
By the construction above,
$\cO_{F^\Inf}$
is stable under the $\cO_{\lsup LG}$-coaction
on $\cO_F$.
In particular, $F^\Inf$ is an $\lsup LG$-torsor
over $\spec \bA_\Inf$.

It remains to show $F^\Inf$ is $\Gal_K$-stable.
It follows from Tannakian formalism and the proof
of \cite[Corollary F.23]{EG23}.
\end{proof}

\subsubsection{Remark}
\label{rem:recover-F-del}
Given the Breuil-Kisin-Fargues $\Gal_K$-lattice
$F^\Inf$,
we can recover the Breuil-Kisin-Fargues $\Gal_E$-lattice
$F^\Inf_\Delta$.
Consider
$$
\bar F^\Inf := F^\Inf \times^{\lsup LG}\Gal(E/K).
$$
We have $\bar F^\Inf \cong \spec \bA\otimes_{\fS_K}\fS_E$,
and such an identification is uniquely determined by
Equation (\ref{eq:wh-G-lvl}).
The following Cartesian diagram
$$
\xymatrix{
F_\Delta^\Inf \ar[d]\ar@{^{(}->}[r] & 
F^\Inf\ar[d]
\\
\spec \bA_\Inf \ar@{^{(}->}[r]^{\hspace{6mm}\Delta} &
\bar F^\Inf
}
$$
recovers $F_\Delta^\Inf$
from $F^\Inf$.

\subsubsection{Definition}
For each Galois extension $L/F$,
define $\cC_{\BKF-L, \lsup LG}^a$
to be the limit-preserving stack over $\cO/\varpi^a$
such that $\cC_{\BKF-L, \lsup LG}^a(A)$
is the groupoid of 
pairs $(F^\Inf, c)$
where $F^\Inf$ is a
Breuil-Kisin-Fargues
$\Gal_L$-module with $\lsup LG$-structure
with $A$-coefficients 
and
$c$ is an $\Gal_L$-equivariant morphism
$F^\Inf\times^{\lsup LG}\Gal(E/K)\cong \spec \bA_{\Inf, A}\otimes_{\fS_K}\fS_E$,
for all finite type
$\cO/\varpi^a$-algebras $A$.

If $L\supset E$, 
there exists a canonical morphism
$$
\cC_{\BKF-K, \lsup LG}^a
\to 
\cC_{\BKF-L, \wh G}^a
$$
defined by sending
$(F^\Inf, c)$ to the pullback of $F^\Inf$
along $c^{-1}\circ\Delta:\spec \bA_{\Inf,A} \to F\times^{\lsup LG}\Gal(E/K)$.

Set
$$
\cC^{L/K,a}_{\lsup L G, \sS, h}:=
\cC_{\BKF-K, \lsup L G}^a
\underset{\cC_{\BKF-L, \wh G}^a}{\times}
\cC_{L, \wh G, \sS, h}^a
$$
and
$$
\cC^{L/K}_{\lsup L G, \sS, h}
:=\dirlim{a}~\cC^{L/K,a}_{\lsup L G, \sS, h}
$$
for all Galois extensions $L/E$.

\subsubsection{Lemma}
\label{lem:BKFL-1}
The morphism
$$
\cC^{L/K,a}_{\lsup L G, \sS, h}
\hookrightarrow
\cX_{K, \lsup LG}
\underset{\cX_{L, \wh G}}{\times}
\cC_{L, \wh G, \sS, h}^a
$$
is a closed immersion.

\begin{proof}
By Remark \ref{rem:recover-F-del}
and the construction before Proposition \ref{prop:ss-BKF},
we see the morphism is a monomorphism.
$\cC^{L/K,a}_{\lsup L G, \sS, h}$
is the $\Gal_K$-stable locus of 
$\cX_{K, \lsup LG}
\underset{\cX_{L, \wh G}}{\times}
\cC_{L, \wh G, \sS, h}^a$, which is a closed condition.
\end{proof}

\subsubsection{Lemma}
For any Galois extension $L/E$,
$\cC^{L/K}_{\lsup L G, \sS, h}$
is a $p$-adic formal algebraic stack
of finite presentation
with affine diagonal.
The forgetful morphism
$\cC^{L/K}_{\lsup L G, \sS, h}\to\cX_{K, \lsup LG}$
is representable by algebraic spaces,
proper and of finite presentation.

\begin{proof}
Combine Lemma \ref{lem:BKF-2}
and Lemma \ref{lem:BKFL-1}.
\end{proof}

\subsection{Inertial types}
\label{subsec:inertial}

Let $A^\circ$ be a $p$-adically complete flat $\cO$-algebra
which is topologically of finite type over $\cO$,
and write $A:=A^\circ[1/p]$.

Let $L/K$ be a finite Galois extension
containing $E$ with inertia group $I_{L/K}$
and suppose $\cO[1/p]$
contains the image of all embeddings $L\hookrightarrow \overline{\bQ}_p$
and that all irreducible $\cO[1/p]$-representations
of $I_{L/K}$ are absolutely irreducible.
Write $l$ for the residue field of $L$
and write $L_0=W(l)[1/p]$.

In \cite[Section 4.6]{EG23},
inertial types 
$\WD_\sigma(\fM^\Inf)$
are attached to Breuil-Kisin-Fargues $\Gal_K$
modules $\fM^\Inf$ with $A^\circ$-coefficients
that admits all descents over $L$.
The underlying $A$-module
of 
$\WD_\sigma(\fM^\Inf)$
is $e_\sigma(\overline{\fM}_{A^\circ} \otimes_{A^\circ}A)$
and it is equipped with an $A$-linear action of $I_{L/K}$
(see loc. cit. for unfamiliar notations).
Here $e_\sigma\in L_0\otimes_{\bQ_p}\cO[1/p]$
is the idempotent corresponding to
a fixed choice of embedding
$\sigma:L_0\hookrightarrow \cO[1/p]$.

\subsubsection{Lemma and Definition}
\label{lem:inertial-1}
Let $A^\circ$ be a $p$-adically complete flat $\cO$-algebra
which is topologically of finite type over $\cO$,
and write $A:=A^\circ[1/p]$.

Let $F^\Inf$ be a Breuil-Kisin-Fargues $\Gal_K$
modules $\fM^\Inf$ with $A^\circ$-coefficients
and $\lsup LG$-structure
that admits all descents over $L$.
The functor $\WD(F^\Inf)$
\begin{align*}
\frep_{\lsup LG} &\to \Vect_A \\
x &\mapsto \WD(F^\Inf \times^{\lsup LG}x)
\end{align*}
is a faithful, exact, symmetric monoidal functor.

\begin{proof}
The functor $\WD(F^\Inf)$ is clearly lax monoidal.
Strict monoidality, faithfulness and exactness are local properties.
Since $A^\circ$ is $\cO$-flat and topologically
of finite type over $\cO$,
it suffices to check the exactness (and monoidality and faithfulness) of $\WD(F^\Inf)$ at $\Lambda$-points
for finite flat $\cO$-algebras $\Lambda$
(because the union of images of $\Lambda$-points
of $\spec A^\circ$ covers all closed points of $\spec A^\circ$).
As is observed in \cite[Remark 4.6.2]{EG23},
in the finite $\cO$-flat coefficients situation,
$\WD(-)$
can be identified with Fontaine's
$D_{\text{pst}}$ functor,
which is well-known to be exact
(see also \cite[Section F.24]{EG23});
faithfulness and monoidality are clear.
\end{proof}

\subsubsection{Lemma}
\label{lem:inertial-2}
Let $A^\circ$ be a $p$-adically complete flat $\cO$-algebra
which is topologically of finite type over $\cO$,
and write $A:=A^\circ[1/p]$.

Let $F^\Inf$ be a Breuil-Kisin-Fargues $\Gal_K$
modules $\fM^\Inf$ with $A^\circ$-coefficients
and $\lsup LG$-structure
that admits all descents over $L$.
Then $\WD(F^\Inf)$
is a $\lsup LG$-torsor over $\spec A$
equipped with an action of $I_F$,
whose formation is compatible with
base change $A^\circ\to B^\circ$
of $p$-adically complete $\cO$-flat algebras
which are topologically of finite type over $\cO$.

\begin{proof}
Combine \cite[Proposition 4.6.3]{EG23}
and Lemma \ref{lem:inertial-1}.
\end{proof}

\subsubsection{Definition}
Let $\tau$ be an $\lsup LG(\cO[1/p])$-valued representation of $I_{L/K}$.
In the setting of 
Lemma \ref{lem:inertial-2},
we say $F^\Inf$
has inertial type $\tau$
if \'etale locally on $\spec A$,
$\WD(F^\Inf)$
is isomorphic to the base change to $A$ of $\tau$.

\subsubsection{Lemma}
\label{lem:fin-num-types}
Fix a finite Galois extension $L/E$.
There are finitely many $\wh G$-conjugacy
classes of inertial types $I_{L/K}\to\lsup LG(\bar\bQ_p)$.

\begin{proof}
Since $I_{L/K}$ is a finite group, and $\bar\bQ_p$ is a characteristic $0$ field,
all group homomorphisms $I_{L/K}\to \lsup LG(\bar\bQ_p[\varepsilon]/\varepsilon^2)$
factor through $\lsup LG(\bar\bQ_p)$:
embed $\lsup LG$ in $\GL_N$ for some $N$,
if $g\in I_{L/K}$ is sent to $x + \varepsilon y$,
then $(x+\varepsilon y)^n=x^n=1$ for some positive integer $n$, and thus $y=0$.
Therefore, deformation theory is trivial in our context.

By the proof of Lemma \ref{lem:inertial-dec} below
(which makes use of only affine GIT theory and there is no circular reasoning),
the coarse moduli space of all $\wh G$-conjugacy
classes of inertial types $I_{L/K}\to\lsup LG(\bar\bQ_p)$
is a finite type scheme over $\bar\bQ_p$.
So we are done.
\end{proof}

\subsubsection{Lemma}
\label{lem:inertial-dec}
In the setting of 
Lemma \ref{lem:inertial-2},
we can decompose $\spec A$
as the disjoint union of open and closed
subschemes $\spec A^\tau$,
where $\spec A^\tau$
is the locus over which $F^\Inf$
has inertial type $\tau$.
Furthermore, the formation of this decomposition
is compatible with
base change $A^\circ\to B^\circ$
of $p$-adically complete $\cO$-flat algebras
which are topologically of finite type over $\cO$.

\begin{proof}
We start with analyzing the (coarse) moduli of inertial types.
Say $I_{L/K}=\{x_1, \dots, x_N\}$ consists of $N$ elements.
Consider the conjugation action of $\wh G_{\bar\bQ_p}$ on
the $N$-tuple
$$
\lsup LG_{\bar\bQ_p}^{(N)}:=\lsup LG_{\bar\bQ_p} \times \cdots \times \lsup LG_{\bar\bQ_p}
$$
Write $\lsup LG_{\bar\bQ_p}^{(N)}\sslash\wh G_{\bar\bQ_p}$
for the GIT quotient.
By \cite[Theorem 1.1]{BMRT11} (which allows disconnected groups),
if $x_1, \dots, x_N$ generate a finite subgroup of $\lsup LG(\bar\bQ_p)$, then
the orbit $\wh G_{\bar\bQ_p}\cdot (x_1, \dots, x_N)$
is closed in $\lsup LG_{\bar\bQ_p}^{(N)}$.
Since the affine GIT quotient is a {\it good quotient},
the image of $\wh G_{\bar\bQ_p}\cdot(x_1, \dots, x_N)$ in the GIT quotient is a closed point.
Inertial types $I_{L/K}\to \lsup LG(\bar \bQ_p)$
corresponds to tuples $(x_1, \dots, x_N)$
that satisfy a finite number of equations
imposed by the group laws;
write $X\subset \lsup LG_{\bar\bQ_p}^{(N)}$
for the closed affine subscheme corresponding to
inertial types.
For each $(x_1, \dots, x_N)\in X$, $\{x_1, \dots, x_n\}$ generate a finite subgroup, and thus $X\sslash\wh G_{\bar\bQ_p}$ is an orbit space.
Since there are only finitely many conjugacy classes
of inertial types by Lemma \ref{lem:fin-num-types},
$X\sslash\wh G_{\bar \bQ_p}$
is zero-dimensional and is thus a disjoint union of points.
By descent, $X\sslash \wh G$ defined over $\cO[1/p]$
is also a zero-dimensional affine variety.
Since the formation of $WD(F^\Inf)$ is compatible
with base change,
the functor of points interpretation
yields a canonical morphism
$\spec A \to [X/\wh G] \to X\sslash\wh G$.
The decomposition $\spec A = \amalg \spec A^\tau$
is the base change of the corresponding decomposition
on $X\sslash\wh G$.
\end{proof}

\subsection{$p$-adic Hodge types}

Next we analyze $p$-adic Hodge types.
By the geometric Shapiro's lemma \cite[Proposition 7.2.4]{L23B},
there is no difference in working with
$G$ or $\Res_{K/\Qp}G$.
For ease of notation, we will often replace $G$ by $\Res_{K/\Qp}G$
and insist $K=\Qp$.

Recall the following characterization of cocharacters.

\subsubsection{Lemma} (\cite[Lemma 3.0.10]{Ba12})
\label{lem:Hodge-tan}
Let $H$ be a split connected reductive group over $E$
and let $\mu,~\mu'$ be two cocharacters of $H_{\bar\bQ_p}$
defined over $E$.
If for any algebraic representation $x:H \to \GL(V)$,
cocharatcers $x\circ \mu$ and $x\circ\mu'$
are conjugate by an element of $\GL(V)(\bar\bQ_p)$,
then $\mu,~\mu'$ are conjugate by an element of $H(\bar\bQ_p)$.

\subsubsection{Corollary}
\label{cor:type-tan}
Let $H$ be a split connected reductive group over $E$.
Let $F$ be a trivial $H$-torsor over $E$.
Let $\eta,~\eta'$ be two exact $\otimes$-filtrations
on $F$.
If for any algebraic representation $x:H \to \GL(V)$,
cocharatcers $x\circ \eta$ and $x\circ\eta'$
are conjugate by an element of $\GL(V)(\bar\bQ_p)$,
then $\eta,~\eta'$ are conjugate by an element of $H(\bar\bQ_p)$.

\begin{proof}
By \cite[IV.2.4]{SN72},
both $\eta, ~\eta'$ are splittable exact $\otimes$-filtrations,
that is, they both are the canonical filtrations attached
to exact $\otimes$-gradings $\wt{\eta}$, $\wt{\eta}'$ on $F$.
An exact $\otimes$-grading $\wt{\eta}$ is equivalent to
a cocharacter $\mu:\Gm\to \underline{\Aut}^{\otimes}(F)$.
Since $F$ is a trivial torsor,
a choice of trivialization induces a group scheme isomorphism
$\underline{\Aut}^{\otimes}(F) \cong H$
and the choice of trivialization does not
affect the conjugacy class of the composition
$\Gm\xrightarrow{\wt{\eta}}\underline{\Aut}^{\otimes}(F) \cong H$, which we also denote by $\wt{\eta}$ by abuse of notation.

We remark that two filtrations $\eta,~\eta'$
are $H$-conjugate if and only if two corresponding cocharacters 
$\mu,~\mu'$ are $H$-conjugate.
Indeed, two cocharacters $\mu_1,~\mu_2$
induce the same filtration if and only if
they are conjugate by
an element of $\underline{\Aut}^{\otimes !}_{\mu}(F)$
which is a closed subgroup scheme of $\underline{\Aut}^{\otimes}(F)$.
See \cite[Section 2.7]{BG19} for unfamiliar notations.

The corollary now follows from Lemma 
\ref{lem:Hodge-tan}.
\end{proof}

\subsubsection{Conjugacy classes of filtrations and Hodge types}
Recall that in \cite[Definition 4.7.7]{EG23},
A {\it Hodge type of rank $d$} $\underline{\lambda}$
is defined to be a set of tuples
of integers $\{\lambda_{\sigma,j}\}_{\sigma:K\hookrightarrow \bar\bQ_p, 1\le j\le d}$
with $\lambda_{\sigma, j} \ge \lambda_{\sigma, j+1}$.
It is clear that rank-$d$ Hodge types $\underline{\lambda}$
are in one-to-one correspondence with
$\wh{\Res_{K/\Qp}\GL_d}$-conjugacy classes
of cocharacters of $\wh{\Res_{K/\Qp}\GL_d}$.

So Hodge types for $G$ should be defined as conjugacy classes of cocharacters
of $\wh{\Res_{K/\Qp}G}$.
As is explained at the beginning of this subsection,
it is harmless to replace $G$ by $\Res_{K/\Qp}G$,
so Hodge types become conjugacy classes of cocharacters
of $\wh{G}$.

\subsubsection{Definition}
\label{def:Hodge}
Let $A^\circ$ be a $p$-adically complete flat $\cO$-algebra
which is topologically of finite type over $\cO$,
and let $F^\Inf$ be a Breuil-Kisin-Fargues $\Gal_E$-module
with $\wh G$-structure of height at most $h$, 
which admits all descents over some finite extension $L$
of $E$.

Let $\underline{\lambda}:\Gm\to \wh G_{\bQp}$
of a cocharacter of $\wh G_{\bQp}$.
We say $F^\Inf$ has {\it Hodge type $\underline{\lambda}$}
if for all algebraic representations $f: \wh G \to \GL(V)$,
$F^\Inf \times^{\wh G} V$
has Hodge type $\underline{\lambda}\circ f$
in the sense of \cite[Corollary 4.7.8]{EG23}.

\subsubsection{Lemma}
\label{lem:Hodge-bd}
Fix a number $C>0$.
Fix an embedding $i:\wh G\hookrightarrow \GL_N$.
There exists only finitely many conjugacy classes
of cocharacters $\underline{\lambda}$ of $\wh G_{\bar\bQ_p}$
such that  $\underline{\lambda}\circ i$
correspond to a tuple of integers $\{\lambda_{i, 1}, \dots, \lambda_{i, N}\}$ such that each $\lambda_{i,\bullet}$ has absolute value bounded by $C$.

\begin{proof}
Clear.
\end{proof}

\subsubsection{Lemma}
\label{lem:Hodge-dec}
In the context of Definition \ref{def:Hodge},
we can write $\spec A$
as a disjoint union of open and closed substacks
$\spec A^{\underline{\lambda}}$
over which $F^\Inf$ is of Hodge type $\underline{\lambda}$.
Moreover, this decomposition is compatible
with base change $A^\circ\to B^\circ$
of $p$-adically complete flat $\cO$-algebras
which are topologically of finite type over $\cO$.

\begin{proof}
By Corollary \ref{cor:type-tan},
for any two non-equivalent Hodge types $\underline{\lambda}$
and $\underline{\lambda}'$,
there exists an algebraic representation $f$
such that $\underline{\lambda}\circ f$
and $\underline{\lambda}\circ f$ are not conjugate.

Since $\spec A$ is topologically of finite type, there are only finitely many Hodge types occuring on $\spec A$
(by Lemma \ref{lem:Hodge-bd}, the claim reduces to the $\GL_N$-case, which is well-known).
In particular, we can choose finitely many different algebraic representations $f_1$, \dots, $f_m$
that distinguish all Hodge types occuring on $\spec A$,
in the sense that for any two non-equivalent Hodge types $\underline{\lambda}$
and $\underline{\lambda}'$,
there exists an algebraic representation $f_i$ ($1\le i\le m$)
such that $\underline{\lambda}\circ f_i$
and $\underline{\lambda}\circ f_i$ are not conjugate.

The lemma has thus been reduced to the $\GL_n$-case,
which is dealt with in \cite[Corollary 4.7.8]{EG23},
since $\spec A^{\underline{\lambda}}$
is the intersection of 
$\spec A^{\underline{\lambda}\circ f_i}$, $1\le i\le m$.
\end{proof}

\subsubsection{Corollary}
In the context of Definition \ref{def:Hodge},
we can write 
$$
\spec A = \coprod_{\underline{\lambda}, \tau}\spec A^{\underline{\lambda}, \tau}
$$
where $\spec A^{\underline{\lambda}, \tau}$ is the locus
over which $F^\Inf$ is of Hodge type $\underline{\lambda}$ and inertial type $\tau$.
Moreover, this decomposition is compatible
with base change $A^\circ\to B^\circ$
of $p$-adically complete flat $\cO$-algebras
which are topologically of finite type over $\cO$.

\begin{proof}
Combine Lemma \ref{lem:Hodge-dec} and Lemma \ref{lem:inertial-dec}.
\end{proof}

\subsubsection{Remark}
(1) By Lemma \ref{lem:Hodge-tan},
there is an intrinsic characterization of Hodge types
for Breuil-Kisin-Fargues modules with $\wh G$-structure.
The Hodge filtration on the de Rham periods
(\cite[Definition 4.7.6]{EG23})
is exact and $\otimes$-compatible,
and thus defines an exact $\otimes$-filtration
$\mu$
on $D_{\operatorname{dR}}(F^\Inf)$.
The clopen $\spec A^{\underline{\lambda}}$
in Lemma \ref{lem:Hodge-dec}
is precisely the locus of $\spec A$ where
$\mu$ is conjugate to $\underline{\lambda}$
at all closed points.

(2)
The formation of Hodge types
is not sensitive to restriction of fields.
So we say a Breuil-Kisin-Fargues $\Gal_K$
module with $\lsup LG$-structure
has Hodge type $\underline{\lambda}$
if when restricted to $\Gal_E$,
    it is a Breuil-Kisin-Fargues $\Gal_E$
module with $\wh G$-structure
and Hodge type $\underline{\lambda}$.

\subsubsection{Definition}
Write $\cC^{L/K, \fl}_{\lsup LG, \sS, h}$
for the maximal $\cO$-flat substack
of $\cC^{L/K}_{\lsup LG, \sS, h}$
(see \cite[Appendix A]{EG23} for the existence).
We record the following proposition.

\subsubsection{Proposition}
\label{prop:Css}
Let $L/E$ be a finite Galois extension.
Then $\cC^{L/K, \fl}_{\lsup LG, \sS, h}$
is the scheme-theoretical union
of closed substacks
$\cC^{L/K, \fl}_{\lsup LG, \sS, h, \underline{\lambda}, \tau}$,
which is uniquely characterized by the following property:
if $A^\circ$ is a finite flat $\cO$-algebra,
    then an $A^\circ$-point of $\cC^{L/K, \fl}_{\lsup LG, \sS, h}$
    factors through $\cC^{L/K, \fl}_{\lsup LG, \sS, h, \underline{\lambda}, \tau}$
    if and only if the corresponding
        Breuil-Kisin-Fargues $F^\Inf$
        has Hodge type $\underline{\lambda}$
        and inertial type $\tau$.

Moreover, $\cC^{L/K, \fl}_{\lsup LG, \sS, h, \underline{\lambda}, \tau}$
is a $p$-adic formal algebraic stack of finite presentation
which is flat over $\spf \cO$ and whose diagonal is affine.
The natural morphism 
$\cC^{L/K, \fl}_{\lsup LG, \sS, h, \underline{\lambda}, \tau}\to \cX_{K, \lsup LG}$
is representable by algebraic spaces, proper, and of finite presentation.

\begin{proof}
The proof is formally identical to that of \cite[Proposition 4.8.2]{EG23}.
\end{proof}

\subsubsection{Definition}
Let $\tau$ be an inertial type
and let $\underline{\lambda}$ be a Hodge type.
Let $h$ be a sufficiently large integer
such that $\underline{\lambda}\circ i$
is bounded by $h$ (where $i$ is the fixed embedding $\wh G\to \GL_d$),
and
let $L/E$ be a Galois extension
such that $\tau$ is trivial when restricted to $I_L$.
Define $\cX_{K, \lsup LG}^{\lambda, \tau}$
to be the scheme-theoretic image
of $\cC^{L/K, \fl}_{\lsup LG, \sS, h, \underline{\lambda}, \tau}$ in $\cX_{K, \lsup LG}$.

\subsection{Potentially semistable deformation rings}~
Recall the following theorem about Galois deformation rings.

\subsubsection{Theorem} (\cite[Theorem A]{BG19}, \cite[Theorem 3.0.12]{Ba12})
\label{thm:def}
Let $\tau$ be an inertial type, and let $\underline{\lambda}$ be a Hodge type.
Fix a mod $p$ $L$-parameter $\bar\rho:\Gal_K\to \lsup LG(\bF)$.
The framed potentially semistable deformation ring
$R_{\bar\rho}^{\square, \tau, \underline{\lambda}}$
is equidimensional of dimension
$$
1+ \dim_{\bQ_p} \wh G + \dim_{\bQ_p} \wh G/P_{\underline{\lambda}}
$$
where $P_{\underline{\lambda}}:=\underline{\Aut}^{\otimes}(\underline{\lambda})$ is the parabolic subgroup of $\wh G$ which stabilizes the exact $\otimes$-filtration associated to $v$.
(See \cite[Section 2.7]{BG19} for unfamiliar notations.)

Denote by $R_{\bar\rho}^\square$ the universal deformation ring.
The potentially semistable deformation ring 
$R_{\bar\rho}^{\square, \tau, \underline{\lambda}}$
is the unique $\cO$-flat quotient
such that for any finite local $\cO[1/p]$-algebra
$B$, any $\cO[1/p]$-homomorphism
$\zeta: R_{\bar\rho}^\square \to B$ factors through
$R_{\bar\rho}^{\square, \tau, \underline{\lambda}}$
if and only if $\zeta$ corresponds to an $L$-parameter
    $\rho_\zeta: \Gal_K\to \lsup LG(B)$
    that is potentially semi-stable of inertial type $\tau$ and Hodge type $\underline{\lambda}$.

\begin{proof}
The first paragraph is \cite[Theorem A]{BG19};
the second paragraph is \cite[Theorem 3.0.12]{Ba12}.
We remark that in \cite{BG19}, any non-split group $G$
is allowed.
In \cite{Ba12}, only split groups $G$ are considered.
However, the same argument works through as long
as deformation rings with fixed inertial types can be constructed,
which we have done in Subsection \ref{subsec:inertial}.
We also remark that \cite{BG19} directly uses results
of \cite{Ba12}
even though the setting of \cite{BG19} is more general.
\end{proof}

\subsubsection{Lemma}
\label{lem:def-fact}
In the setting of Theorem \ref{thm:def},
the morphism $\spf R_{\bar\rho}^{\square, \tau, \underline{\lambda}}\to \cX_{K, \lsup LG}$
factors through
a versal morphism
$$
\spf R_{\bar\rho}^{\square, \tau, \underline{\lambda}}
\to \cX^{\tau, \underline{\lambda}}_{K, \lsup LG}.
$$

\begin{proof}
The proof is formally identical to that
of \cite[Proposition 4.8.10]{EG23}.
The inputs are
\begin{itemize}
\item Algebraicity of $\cX_{K, \lsup LG}$ (the main theorem of \cite{L23B}, enhanced by \cite{Min24});
\item Description of finite $\cO$-flat points of $\cC^{L/K, \fl}_{\lsup LG, \sS, h, \underline{\lambda}, \tau}$ (Proposition \ref{prop:Css});
\item Description of finite $\cO$-flat points of $\spf R_{\bar\rho}^{\square, \tau, \underline{\lambda}}$ (Theorem \ref{thm:def}); and
\item Existence of Breuil-Kisin-Fargues lattices.
\end{itemize}
The last bullet point in the $\GL_d$-case is \cite[Theorem 4.7.13]{EG23},
which holds for general (split) groups $G$, as is observed in
the thesis of B. Levin (see the proof of \cite[Proposition 5.4.2]{Lev13}).
The non-split group case follows from the split group case
(see the proof of Proposition  \ref{prop:ss-BKF}).
\end{proof}

\subsubsection{Theorem}
\label{thm:ss}
Let $\tau$ be an inertial type
and let $\underline{\lambda}$ be a Hodge type.
Then $\cX_{K, \lsup LG}^{\tau, \underline{\lambda}}$
is a $p$-adic formal algebraic stack
which is of finite type and flat over $\spf \cO$.
It is uniquely determined as
the $\cO$-flat closed substack of $\cX_{K, \lsup LG}$
by the following property:
if $A^\circ$ is a finite $\cO$-flat algebra,
then  $\cX_{K, \lsup LG}^{\tau, \underline{\lambda}}(A^\circ)$
is the subgroupoid consisting of $L$-parameters
which become potentially semistable
of Hodge type $\underline{\lambda}$
and inertia type $\tau$ after inverting $p$.

The mod $p$ fiber
$$
\cX_{K, \lsup LG}^{\tau, \underline{\lambda}}
\underset{\spf \cO}{\times}\spec \bF
$$
is equidimensional of dimension 
$\dim_{\bQ_p} \wh G/P_{\underline{\lambda}}.$

\begin{proof}
The proof is formally identical to that of \cite[Theorem 4.8.12]{EG23}.
The first claim follows from Proposition \ref{prop:Css},
\cite[Theorem 2]{L23B} and \cite[Proposition A.21]{L23B}.
The description of finite $\cO$-flat points
$\cX_{K, \lsup LG}^{\tau, \underline{\lambda}}(A^\circ)$
follows from Lemma \ref{lem:def-fact}.
The uniqueness follows from \cite[Proposition 4.8.6]{EG23}.
The dimension calculation follows from Theorem \ref{thm:def}
and the versality of 
$\spf R_{\bar\rho}^{\square, \tau, \underline{\lambda}}$
established in Lemma \ref{lem:def-fact} (compare with \cite[Theorem 4.8.14]{EG23}).
\end{proof}

\end{appendices}

\addcontentsline{toc}{section}{References}

\bibliographystyle{amsalpha}
\bibliography{unobs}

\end{document}